\documentclass[12pt,reqno]{amsart}
\usepackage{amssymb}
\input amssym.def
\usepackage{amsmath,amsfonts,hyperref,xcolor,mathtools}
\usepackage{amscd,accents}
\usepackage[mathscr]{eucal}
\usepackage{enumitem}
\usepackage{orcidlink}
\usepackage[utf8]{inputenc}

\hypersetup{colorlinks=true,citecolor={purple},linkcolor={teal},urlcolor={violet}}

\newcommand{\sumdash}{\mathop{{\sum}'}}
\newcommand{\E}{{\mathbb E}}
\newcommand{\N}{{\mathbb N}}

\newcommand{\C}{{\mathbb C}}
\newcommand{\R}{{\mathbb R}}

\newtheorem{thm}{Theorem}[section]
\newtheorem{lem}{Lemma}[section]

\newtheorem{rmk}{Remark}[section]

\newtheorem{conj}{Conjecture}[section]

\newcommand{\thmref}[1]{Theorem~\ref{#1}}

\newcommand{\lemref}[1]{Lemma~\ref{#1}}

\makeatletter
\@namedef{subjclassname@2020}{%
\textup{2020} Mathematics Subject Classification}
\makeatother

\begin{document}

\date{\today}

\title[Triple convolution sums of the generalised divisor functions]
{Triple convolution sums of the generalised divisor functions and related
sums over primes}

\author{Bikram Misra \orcidlink{0009-0000-5863-2789}
and Biswajyoti Saha \orcidlink{0009-0009-2904-4860}}

\address{Bikram Misra and Biswajyoti Saha\ \newline
Department of Mathematics, Indian Institute of Technology Delhi, 
New Delhi 110016, India.}
\email{bikram.misra@maths.iitd.ac.in, biswajyoti@maths.iitd.ac.in}



\begin{abstract}
We study the triple convolution sum of the generalised 
divisor functions
$$\sum_{n\leq x} d_k(n+h)d_l(n)d_m(n-h),$$
where $h \le x^{1-\epsilon}$ for any $\epsilon>0$ and $d_k(n)$ denotes the 
generalised divisor function which counts the number of ways $n$ can
be written as a product of $k$ many positive integers.
The purpose of this paper is three-fold.  Firstly, we 
note a predicted asymptotic estimate for the above sum,
where the constant appearing in the estimate
can be obtained from the theory of Dirichlet series of 
several complex variables and also using some probabilistic arguments.
Then we show that a lower bound of the correct order can be derived
using the several variable Tauberian theorems, where, more
importantly, the constant in the predicted asymptotic can be recovered.
Lastly, in the spirit of the Titchmarsh divisor problem, we consider this
triple convolution sum over the prime numbers, which essentially leads to a
shifted convolution sum. We use the Tauberian theory of multiple Dirichlet series along with the Bombieri-Vinogradov theorem to derive an explicit lower bound of
this.
\end{abstract}

\maketitle

\section{Introduction}
Throughout this article, $\N$ denotes the set of positive integers.
The study of the convolution sums of arithmetic functions is a fundamental topic of research.
In fact, the study of the convolution sums of the (generalised) divisor functions itself has a towering history,
starting with Ingham's work on the shifted and the additive convolution sums of the divisor function.
For $n\in \N$, let $d(n)$ denote the number of positive divisors of $n$.
Ingham \cite{AEI} showed that for a positive integer $h$,
\begin{equation}\label{shifted-divisor}
\sum_{n \le N} d(n) d(n+h) = \frac{6}{\pi^2} \sigma_{-1}(h) N (\log N)^2+ O(N \log N),
\end{equation}
as $N \to \infty$ and 
\begin{equation}\label{add-divisor-1}
\sum_{n < N} d(n) d(N-n) = \frac{6}{\pi^2} \sigma_{1}(N) (\log N)^2 + O(\sigma_1(N) \log N \log\log N),
\end{equation}
as $N \to \infty$, where $\sigma_s(n):=\sum_{d \mid n\atop d>0} d^s$ for a complex number $s$.
He then used \eqref{shifted-divisor} to study the fourth moment of the Riemann zeta function on the line $\Re(s)=1/2$.
The formulae \eqref{shifted-divisor} and \eqref{add-divisor-1} have subsequently been extended in many directions.
We focus our attention to the ones for the (higher) convolution sums of the (generalised) divisor functions. For $n, k \in \N$, let
$$
d_k(n) = \#\left\{(x_1,\hdots,x_k)\in \N^k : x_1\cdots x_k=n\right\}.
$$
Hence for $k=2$, we get back the divisor function. Clearly, $d_k(n)$ is the coefficient of $n^{-s}$ in the Dirichlet series of $\zeta(s)^k$.
The problem of finding an asymptotic formula for the shifted convolution of $d_k$ is still open for $k\ge 3$.
In 2001, Conrey and Gonek \cite{Gonek-Conrey} gave a conjectural asymptotic formula of the sum $\sum_{n\le x}d_k(n)d_k(n+h)$.
In his blog, Tao \cite{Tao} predicted the following generalised version of this conjecture, using probabilistic arguments.

\begin{conj}\label{tao's conjecture}
Let $\epsilon>0$ and integers $k,l\ge 2$. For $1\le h\le x^{1-\epsilon}$, we have 
\begin{equation}\label{tao's conjecture statement}
\sum_{n\le x}d_k(n+h)d_l(n) \sim \frac{c_{k,l}(h)}{(k-1)!(l-1)!}x(\log x)^{k+l-2}
\end{equation}
as $x\to\infty$, where $c_{k,l}(h)$ is an explicit constant.
\end{conj}

While the conjecture is still open, Ng and Thom \cite{Ng} proved an explicit lower bound of the correct order. 
In \cite{Ng}, they mentioned that Daniel \cite{Daniel} and Henriot \cite{H} in their unpublished works
had proved an explicit upper bound of the appropriate order.
%

The study of the higher convolution of the (generalised) divisor functions is more challenging.
Unlike the shifted convolution sum of the divisor function, an asymptotic estimate for the
triple convolution sum of the divisor function is yet unproven.
For a fixed positive integer $h$, the triple convolution sum of the divisor function is defined as
$$
\mathcal{T}(d,d,d;x,h) := \sum_{n\le x}d(n+h)d(n)d(n-h).
$$
In \cite{TB}, Browning suggested that $\mathcal{T}(d,d,d,x;h) \sim c_h x(\log x)^3$ as $x\to\infty$, with a precise constant $c_h>0$.
This conjecture is still unproven, but Browning \cite{TB} proved an `average' version of the above asymptotic,
which was later extended by Blomer \cite{VB} using spectral tools. A far reaching generalisation of these
average results for the higher convolution sums of a fixed generalised divisor function has recently been obtained by
Matom\"aki, Radziwi\l\l, Shao, Tao and Ter\"av\"ainen \cite{MRSTT}. This also extends a recent work of Miao \cite{miao}.

Using the theory of multiple Dirichlet series, the present authors, along with Murty \cite{MMS}, showed that for a fixed positive integer $h$, as $x \to \infty$,
$$
\mathcal{T}(d,d,d;x,h) \geq c_h x(\log x)^3/27 + O_h(x(\log x)^2).
$$
Note that the constant $c_h$ appearing in the above estimate is exactly the same as in Browning's conjecture.
For the upper bound, it is possible to derive one of the correct order,
from the works of Wolke \cite{DW} and of Nair \cite{MN},
as one has
$$
d(n+h)d(n)d(n-h)\ll_h d((n+h)n(n-h)),
$$
for any $h\neq 0$ (see \cite[Lemma 2.5]{MMS}).

%

In this article, we study the triple convolution sums of the generalised divisor functions, namely,
$$
\mathcal{T}(d_k,d_l,d_m;x,h) := \sum_{n\le x}d_k(n+h)d_l(n)d_m(n-h),
$$
for integers $k,l,m \ge 2$. We first
record an expected asymptotic for $\mathcal{T}(d_k,d_l,d_m;x,h)$,
where the constant appearing in the asymptotic formula, in fact,
can be obtained from the theory of Dirichlet series of 
several complex variables and also using some probabilistic arguments,
as we shall see below.

\begin{conj}\label{conj-nabla}
Let $\epsilon>0$ and integers $k,l,m\ge 2$. Then for $0<h\le x^{1-\epsilon}$,
\begin{equation}\label{our-conjecture-statement}
\mathcal{T}(d_k,d_l,d_m;x,h) \sim \nabla_{k,l,m}(h)\frac{x(\log x)^{k+l+m-3}}{(k-1)!(l-1)!(m-1)!},
\end{equation} 
as $x \to \infty$, where the constant $\nabla_{k,l,m}(h)$ is given as follows:
\begin{align*}
&\nabla_{k,l,m}(h)= C_{k,l,m}(h)D_{k,l,m}(h)\prod_{p}\left(1-\frac{1}{p}\right)^{k+l+m-3},
\end{align*} 
with
\begin{align*}
C_{k,l,m}(h) = \prod_{p\mid 2h}\left(\sum_{\nu_1,\nu_2,\nu_3\ge 0}d_{k-1}(p^{\nu_1})d_{l-1}(p^{\nu_2})d_{m-1}(p^{\nu_3})\frac{g(p^{\nu_1},p^{\nu_2},p^{\nu_3})}{[p^{\nu_1},p^{\nu_2},p^{\nu_3}]}\right),
\end{align*}
and
\begin{align*}
D_{k,l,m}(h)=&\prod_{p\ \nmid \ 2h}\left\{\left(1-\frac{1}{p}\right)^{1-k}+\left(1-\frac{1}{p}\right)^{1-l}+\left(1-\frac{1}{p}\right)^{1-m}-2\right\}.
\end{align*}
Here $g$ is the function such that $g(u,v,w) = 1$ if the system $n\equiv-h\bmod u,n\equiv0\bmod v,n\equiv h\bmod w$, has a solution, else it is 0 and $[u,v,w]$ denotes the least common multiple of $u,v$ and $w$.
\end{conj}

As mentioned, the appearance of the constant in \eqref{our-conjecture-statement} can be motivated
using some probabilistic arguments and by the theory of Dirichlet series of several complex variables. Regarding this expected asymptotic,
we first establish a lower bound of the correct order, with the constant emerging from the theory of
Dirichlet series of several complex variables (also see \thmref{three-div-function}).
The probabilistic heuristic is discussed in the last section.

\begin{thm}\label{lower bound of triple}
Let $\epsilon>0$ and integers $k,l,m\ge 2$. Then for $h\le x^{1-\epsilon}$,
\begin{equation}
\mathcal{T}(d_k,d_l,d_m;x,h) \ge \frac{\nabla_{k,l,m}(h)}{3^{k+l+m-3}}\frac{x(\log x)^{k+l+m-3}}{(k-1)!(l-1)!(m-1)!} + O(x(\log x)^{k+l+m-4}),
\end{equation} 
as $x \to \infty$.
\end{thm}

When we consider $k=l=m$, again it is possible to derive an upper bound of the appropriate order, using
the work of Wolke \cite{DW} and of Nair \cite{MN}, as here also
$$
d_k(n+h)d_k(n)d_k(n-h)\ll_{h,k} d_k((n+h)n(n-h)),
$$
for any $h\neq 0$ and $k \ge 2$ (see the proof of \cite[Lemma 2.5]{MMS}). For an upper bound of the correct order, for general choices of $k,l$ and $m$, we can use the work of Nair and Tenebaum \cite{NT}, which provides a vast generalisation of Nair's work \cite{MN}. Henriot \cite{H-1} extended the work of Nair and Tenebaum, and proved a bound that is uniform with respect to the discriminant.

Now in the case of summing over the prime numbers, the probelm of
	deriving an asymptotic formula for the sum $\sum_{p\le x}d(p-a)$, for a fixed non-zero integer $a$, is known as the Titchmarsh divisor problem. This sum was first studied by Titchmarsh \cite{Titchmarsh}, under the assumption of the generalised Riemann hypothesis. In 1961, Linnik \cite{Linnik} applied the dispersion method to treat the sum $\sum_{p\le x}d(p-a)$. Rodriquez \cite{Rodriquez} and Halberstam \cite{Halberstam} independently proved the asymptotic formula by applying the celebrated Bombieri-Vinogradov Theorem. 
	
	In spirit of the Titchmarsh divisor problem, it is natural to consider the following triple convolution sum over primes,
    $$
    \mathcal{T}'(d_k,d_l,d_m;x,h):=\sum_{p\le x}d_k(p+h)d_l(p)d_m(p-h),
    $$
    for $k,l,m\ge 2$. 
    As $d_l(p)=l$, the triple convolution sum $\mathcal{T}'(d_k,d_l,d_m;x,h)$ essentially leads to the study of the following shifted convolution sum 
    $$
    \mathcal{T}'(d_k,d_m;x,h) :=\sum_{p\le x}d_k(p+h)d_m(p-h).
    $$
    However, the research about the convolution sums over primes has been limited. Assuming a Siegel–Walfisz type theorem on primes in arithmetic progressions with large moduli, Hooley \cite{Hooley-77} proved the asymptotic formula 
    $$
    \sum_{p\le x}d(p+1)d(p+2)\sim \frac{9x\log x}{\pi^2},
    $$
    as $x\to\infty.$
    An upper bound of the appropriate order,
    namely, $\mathcal{T}'(d_k,d_m;x,h) \ll_{h,k,m} x (\log x)^{k+m-3}$, follows from \cite[Theorem 3]{NT}. Using the Tauberian theory of multiple Dirichlet series along-with a variant of the Bombieri-Vinogradov theorem, we derive the following explicit lower bound.
	
	\begin{thm}\label{lower bound of triple sum over primes}
		Let $\epsilon>0$ and integers $k,m\ge 2$. Then for $h\le x^{1-\epsilon}$,
		\begin{equation}
			\mathcal{T}'(d_k,d_m;x,h)\ge \frac{\nabla_{k,m}'(h)}{4^{k+m-2}}\frac{x(\log x)^{k+m-3}}{(k-1)!(m-1)!} + O(x(\log x)^{k+m-4}),
		\end{equation}
		as $x \to \infty$, where
		\begin{equation}\label{constant}
			\nabla_{k,m}'(h)= C_{k,m}'(h)D_{k,m}'(h)\prod_{p}\left(1-\frac{1}{p}\right)^{k+m-2},
		\end{equation} 
		for some arithmetic constants $C_{k,m}'(h)$ and $D_{k,m}'(h)$ defined as follows:
		\begin{align*}
			C_{k,m}'(h)= \prod_{p \ \mid \ 2h}\left(\sum_{\nu_1,\nu_2\ge 0}d_{k-1}(p^{\nu_1})d_{m-1}(p^{\nu_2})\frac{g_1(p^{\nu_1},p^{\nu_2})}{\phi([p^{\nu_1},p^{\nu_2}])}\right)
		\end{align*}
		and
		\begin{align*}
			D_{k,m}'(h)
			=&\prod_{p\ \nmid \ 2h}\left(1-2\left(1-\frac{1}{p}\right)^{-1}+\left(1-\frac{1}{p}\right)^{-k}+\left(1-\frac{1}{p}\right)^{-m}\right).
		\end{align*}
       Here $\phi$ denotes the Euler totient function and $g_1(u,w) = 1$ if the system $n\equiv-h\bmod u,n\equiv h\bmod w$ has a solution, else it is $0$. 
        \end{thm}

    \begin{rmk}\label{rmk}\rm
        The constant $C_{k,m}'(h)$ can also be defined multiplicatively as follows:
        \begin{align*}
         C_{k,m}'(2^{\alpha}) = \ & 1+d_{k}(2^{\alpha+1})\sum_{\nu>\alpha+1}\frac{d_{m}(2^{\nu})}{2^{\nu}}+d_{m}(2^{\alpha+1})\sum_{\nu>\alpha+1}\frac{d_{k}(2^{\nu})}{2^{\nu}}-\frac{d_{m}(2^{\alpha+1})d_{k}(2^{\alpha})}{2^{\alpha}}\\
         & + \sum_{\nu=1}^{\alpha+1}\frac{d_k(2^{\nu})d_m(2^{\nu})-d_k(2^{\nu-1})d_m(2^{\nu-1})}{2^{\nu-1}},
        \end{align*}
        and for an odd prime $p$,
        \begin{align*}
            C_{k,m}'(p^{\alpha}) = \ &  1+d_{m}(p^{\alpha})\sum_{\nu>\alpha}\frac{d_{k}(p^{\nu})}{p^{\nu}}+d_{k}(p^{\alpha})\sum_{\nu>\alpha}\frac{d_{m}(p^{\nu})}{p^{\nu}}-2\left(1-\frac{1}{p}\right)^{-1}\frac{d_{k}(p^{\alpha})d_m(p^{\alpha})}{p^{\alpha+1}}\\
            &+\left(1-\frac{1}{p}\right)^{-1}\sum_{\nu=1}^{\alpha}\frac{d_k(p^{\nu})d_m(p^{\nu})-d_k(p^{\nu-1})d_m(p^{\nu-1})}{p^{\nu}}.
        \end{align*}
	\end{rmk}

\section{Preliminaries}
In this section, we collect the results that are required to prove our theorems.

\subsection{Elementary facts and key results related to the generalised divisor functions}
	We first record some elementary facts and key results related to the generalised divisor functions.
	\begin{itemize}
		\item[(i)] We can write $d_k(n) = \sum_{a|n}d_{k-1}(a)$ for any $k \ge 2$, where $d_1(n)=1$ for all $n \ge 1$.
		\item[(ii)] For a non-negative integer $j$, any divisor of $p^j$ is of the form $p^{x_i}$ for some $0\le x_i\le j$. Therefore $d_k(p^j)$ is same as counting the non-negative integer solutions of the equation
		$$
		x_1+x_2+\cdots+x_k=j.
		$$ 
		In other words,
		\begin{equation}\label{formula of d_k(p^j)}
		d_k(p^j)=\binom{k+j-1}{j}. 
		\end{equation}
		\item[(iii)] For a non-negative integer $j$, we have the identity
		$$
		d_k(p^{j+1})=d_k(p^{j})+d_{k-1}(p^{j+1}).
		$$
		\item[(iv)] For a non-negative integer $j$, applying the above identity twice we get
		\begin{equation}\label{recursive-relation}
        \begin{split}
        & d_k(p^{j+1})d_m(p^{j+1})\\
        & =d_{k-1}(p^{j+1})d_{m-1}(p^{j+1})+ d_k(p^{j})d_m(p^{j})+d_k(p^{j})d_{m-1}(p^{j+1})+ d_m(p^j)d_{k-1}(p^{j+1}).
        \end{split}
		\end{equation}
        \item[(v)] For a non-negative integer $\alpha$, we also have
        \begin{equation}\label{h}
\sum_{j>\alpha}\frac{d_{k-1}(p^{j})}{p^{j}}=
\left(1-\frac{1}{p}\right)\sum_{j>\alpha}\frac{d_{k}(p^{j})}{p^{j}} - \frac{d_{k}(p^{\alpha})}{p^{\alpha+1}}.
\end{equation}
		\item[(vi)] For any prime number $p$ and $k \ge 2$,
        using \eqref{formula of d_k(p^j)}, we have the following power series
        expansion of $(1-x)^{-k}$ for $|x|<1$,
		\begin{equation}\label{Power series expression}	
        \sum_{j=0}^{\infty}d_k(p^j)x^j = (1-x)^{-k}.
		\end{equation}

	\end{itemize}
The following variant of the Chinese remainder theorem can be found in
\cite[Theorem 3-12]{WJL}.  

\begin{lem}\label{CRT}
For positive integers $d_1,\ldots,d_k$ and integers $a_1,\ldots,a_k$, the system 
\begin{equation}\label{lemma-murty-goel}
	    \begin{cases}
		x\equiv a_1\bmod{d_1},\\
		\ \ \ \vdots\\
		x\equiv a_k\bmod{d_k},
	\end{cases}\
\end{equation}
has a solution if and only if $\gcd (d_i,d_j)\mid(a_i-a_j)$ for all $1\le i,j\le k$.
Moreover, when a solution exists, it is unique modulo the least common multiple $[d_1,\ldots,d_k]$ of $d_1,\ldots,d_k$.
\end{lem}

\subsection{Results from the multiple Dirichlet series theory}
This subsection is to recall the Tauberian theorems of de la Bret\`eche \cite{RB}.
Let $\R^+$ denote the set of all non-negative real numbers and $\R_*^+$ denote the set of all positive real numbers.
For a positive integer $m$, we denote an $m$-tuple $(s_1,\ldots,s_m)$ of complex numbers by $\pmb{s}$. Let $\tau_j=\Im(s_j)$ and
$\mathcal{L}_m(\mathbb{C})$ be the space of all linear forms on $\mathbb{C}^m$ over $\mathbb{C}$.
We denote by $\left\{e_j\right\}_{j=1}^{m}$, the canonical basis of $\mathbb{C}^m$ and $\left\{e_j^*\right\}_{j=1}^{m}$,
the dual basis in $\mathcal{L}_m(\mathbb{C})$. Let $\mathcal{L\mathbb{R}}_m(\mathbb{C})$
(respectively $\mathcal{L\mathbb{R}}^+_m(\mathbb{C})$) denote the set of linear forms of $\mathcal{L}_m(\mathbb{C})$ that are having values in $\mathbb{R}$
(respectively, in $\mathbb{R^+}$) when we restrict to
 $\mathbb{R}^m$ (respectively to $(\mathbb{R^+})^m$). Let $\beta_j>0$ for all $j=1,\ldots,m$. Then we denote by $\mathcal{B}$,
the linear form $\sum_{j=1}^{m}\beta_je_j^*$, and $\pmb{\beta}=(\beta_1,\ldots,\beta_m)$ be the associated row matrix.
We define $X^{\pmb{\beta}} := (X^{\beta_1},\ldots,X^{\beta_m})$. Let $\mathcal{L}$ be a family of linear forms and
for this we define $\text{conv}(\mathcal{L}) := \sum_{\ell\in \mathcal{L}}\mathbb{R^+}\ell$ and
$\text{conv}^*(\mathcal{L}) := \sum_{\ell\in \mathcal{L}}\mathbb{(R^*)^+}\ell$.
With these notations in place, \cite[Th\'eor\`eme 1]{RB} reads as follows:

\begin{thm}\label{theorem 1, breteche}
Let $f$ be an arithmetic function on $\mathbb{N}^m$ taking positive values and $F$ be the associated Dirichlet series
$$
F(\pmb{s}) = \sum_{d_1=1}^{\infty}\ldots\sum_{d_m=1}^{\infty}\frac{f(d_1,\ldots,d_m)}{d_1^{s_1}\ldots d_m^{s_m}}. 
$$	
Suppose that there exists $\pmb{\alpha}=(\alpha_1,\ldots, \alpha_m) \in (\mathbb{R^+})^m$ such that $F$ satisfies the following three properties:
\begin{enumerate}
\item The series $F(\pmb{s})$ is absolutely convergent for $\pmb{s} \in \mathbb{C}^m$ such that $\Re(s_i)>{\alpha_i}$.
\item There exists a family $\mathcal{L}$ of $n$ many non-zero linear forms $\mathcal{L}:= \left\{\ell^{(i)}\right\}_{i=1}^{n}$
in $\mathcal{L\mathbb{R}}^+_m(\mathbb{C})$ and a family of finitely many linear forms $\left\{h^{(k)}\right\}_{k\in \mathcal{K}}$
in $\mathcal{L\mathbb{R}}^+_m(\mathbb{C})$, such that the function $H$ from $\mathbb{C}^m$ to $\mathbb{C}$ defined by
$$
H(\pmb{s}) := F(\pmb{s}+\pmb{\alpha})\prod_{i=1}^{n}\ell^{(i)}(\pmb{s}) 
$$
can be extended to a holomorphic function in the domain
$$
\mathcal{D}(\delta_1,\delta_3) := \left\{\pmb{s}\in \mathbb{C}^m: \Re\left(\ell^{(i)}(\pmb{s})\right)>-\delta_1
\ \text{for all} \ i \text{ and } \Re\left(h^{(k)}(\pmb{s})\right)>-\delta_3 \text{ for all } k \in \mathcal K \right\}
$$
for some $\delta_1,\delta_3>0$.
\item There exists $\delta_2>0$ such that for every $\epsilon,\epsilon{'}>0$, the upper bound
$$
|H(\pmb{s})| \ll \left(1+||\Im(\pmb{s})||_{1}^{\epsilon}\right)
\prod_{i=1}^{n} \left(|\Im\left(\ell^{(i)}(\pmb{s})\right)|+1\right)^{1-\delta_2\min\left\{0,\Re\left(\ell^{(i)}(\pmb{s})\right)\right\}}
$$
is uniformly valid in the domain $\mathcal{D}\left(\delta_1-\epsilon{'},\delta_3-\epsilon{'}\right)$.  
\end{enumerate}
Let $J(\pmb{\alpha}) := \left\{j \in \{1,\ldots,m\} : \alpha_j=0\right\}$.
Let $r:= |J(\pmb{\alpha})|$ and $\ell^{(n+1)},\ldots,\ell^{(n+r)}$ be the linear forms $e_j^{*}$ for $j \in J(\pmb{\alpha})$.
Then for ${\pmb \beta}=(\beta_1,\ldots,\beta_m) \in (\mathbb{R^+})^m$,
there exists a polynomial $Q_{\pmb\beta}\in \mathbb{R}[X]$ of degree less than or equal to
$n+r-rank\left(\left\{\ell^{(i)}\right\}_{i=1}^{n+r}\right)$ and a real number
$\theta=\theta\left(\mathcal{L},\left\{h^{(k)}\right\}_{k\in \mathcal{K}},\delta_1,\delta_2,\delta_3,\pmb{\alpha},\pmb{\beta}\right)>0$
such that we have, for $X\ge1$,	
$$
S(X^{\pmb\beta}) := \sum_{1\le d_1\le X^{\beta_1}}\ldots\sum_{1\le d_m\le X^{\beta_m}}f(d_1,\ldots,d_m)
=X^{\langle\pmb{\alpha},\pmb{\beta}\rangle}\left(Q_{\pmb\beta}(\log X)+O(X^{-\theta})\right).
$$
\end{thm}

We also need \cite[Th\'eor\`eme 2]{RB}.

\begin{thm}\label{theorem 2, breteche}
Let the notations be as in \thmref{theorem 1, breteche}. If we have
$\mathcal{B}$ is not in the span of $\left\{\ell^{(i)}\right\}_{i=1}^{n+r}$, then $Q_{\pmb\beta}=0$.
Next suppose thata we have the following two conditions:
\begin{enumerate}
\item There exists a function $G$ such that $H(\pmb{s}) = G(\ell^{(1)}(\pmb{s}),\ldots,\ell^{(n+r)}(\pmb{s}))$.
\item $\mathcal{B}$ is in the span of $\left\{\ell^{(i)}\right\}_{i=1}^{n+r}$ and there exists no strict
subfamily $\mathcal{L'}$ of $\left\{\ell^{(i)}\right\}_{i=1}^{n+r}$ such that $\mathcal{B}$ is in the span of 
$\mathcal{L'}$ with
$$
card(\mathcal{L'})-rank(\mathcal{L'}) = card\left(\left\{\ell^{(i)}\right\}_{i=1}^{n+r}\right)-rank\left(\left\{\ell^{(i)}\right\}_{i=1}^{n+r}\right).
$$
\end{enumerate}
Then, for $X \ge 3$, the polynomial $Q_{\pmb \beta}$ satisfies the relation 
$$
Q_{\pmb\beta}(\log X) = C_0X^{-\langle\pmb{\alpha},\pmb{\beta}\rangle}I(X^{\pmb\beta})+O((\log X)^{\rho-1}),
$$
where $C_0:= H(0,\ldots,0), \rho:= n+r-rank\left(\left\{\ell^{(i)}\right\}_{i=1}^{n+r}\right)$, and
$$
I(X^{\pmb\beta}) := \int\int\ldots\int_{A(X^{\pmb \beta})} \frac{dy_1\ldots dy_n}{\prod_{i=1}^{n} y_i^{1-\ell^{(i)}(\pmb{\alpha})}},
$$ 
with
$$
A(X^{\pmb{\beta}}) := \left\{\pmb{y} \in [1,\infty)^n: \prod_{i=1}^{n}y_i^{\ell^{(i)}(e_j)}\le X^{\beta_j}  \ \text{for all} \  j\right\}.
$$ 
\end{thm}
\subsection{Bombieri-Vinogradov theorem}
We first recall the following notations for the celebrated Bombieri-Vinogradov theorem. For coprime positive integers $a,q$, let $\pi(y;q,a)$ denote the number of primes $p\le y$ such that $p\equiv a\bmod{q}$  and $q\le y$. Let $\phi$ denote the Euler-$\phi$ function and ${\rm li }(y) := \int_{2}^{y} {dt}/{\log t}$. The following version of the Bombieri-Vinogradov theorem is taken from \cite[Theorem 9.2.1]{CM}.

\begin{thm}\label{BV thm}
	For any $A>0$, there exists $B=B(A)>0$ such that as $x \to \infty$,
	\begin{equation}\label{BV}
		\sum_{q\le \frac{x^{1/2}}{(\log x)^B}} \max_{y\le x}\max_{\gcd(a,q)=1}\left|\pi(y;q,a)-\frac{{\rm li }(y)}{\phi(q)}\right| \ll \frac{x}{(\log x)^A}.
	\end{equation}
\end{thm}
    
We also record the following variant of the Bombieri-Vinogradov theorem
that can be found in \cite[Exercise 26.2]{DK} (as $d_k(q)d_m(q) \le d_{km}(q)$).

	\begin{thm}\label{BV-a}
		For integers $k,m\ge 2$ and given $A>0$, there exists $B=B(A,k,m)>0$ such that 
		\begin{equation}
			\sum_{q\le \frac{x^{1/2}}{(\log x)^B}}d_k(q)d_m(q) \max_{y\le x}\max_{\gcd(a,q)=1}\left|\pi(y;q,a)-\frac{{\rm li }(y)}{\phi(q)}\right| \ll \frac{x}{(\log x)^A}.
		\end{equation}
	\end{thm} 

\section{Proof of \thmref{lower bound of triple}}
For the proof of Theorem \ref{lower bound of triple},
we follow the method formulated in \cite{MMS} closely, but we also need to make several key changes. We write the triple convolution sum as 
	\begin{align*}
		\mathcal{T}(d_k,d_l,d_m;x,h) &= \sum_{n\le x}d_k(n+h)d_l(n)d_m(n-h)\\
		&=\sum_{n\le x}\left(\sum_{u\mid n+h}d_{k-1}(u)\right)\left(\sum_{v\mid n}d_{l-1}(v)\right)\left(\sum_{w\mid n-h}d_{m-1}(w)\right)\\
		&=\sum\limits_{\substack{u\le x+h \\ v\le x \\ w\le x-h}}d_{k-1}(u)d_{l-1}(v)d_{m-1}(w)\sumdash_{n\le x}1,
	\end{align*}
where the primed sum denotes the sum for $n \le x$ where $n+h\equiv 0 \bmod u, n\equiv 0\bmod v$ and $n-h\equiv 0\bmod w$. Then
	$$\mathcal{T}(d_k,d_l,d_m;x,h) = S_1(x;h)+S_2(x;h)-S_3(x;h),$$
	where $S_i(x;h)$ are defined as follows:
	\begin{align*}
		&S_1(x;h) := \sum\limits_{{u,v,w\le x}}d_{k-1}(u)d_{l-1}(v)d_{m-1}(w) \ \sumdash_{n\le x}1,\\
		&S_2(x;h) := \sum\limits_{{x<u\le x+h \atop v,w\le x}}d_{k-1}(u)d_{l-1}(v)d_{m-1}(w)\ \sumdash_{n\le x}1,\\
		&S_3(x;h) :=  \sum\limits_{{u\le x+h, v\le x \atop x-h<w\le x}}d_{k-1}(u)d_{l-1}(v)d_{m-1}(w)\ \sumdash_{n\le x}1.
	\end{align*}
	In order to derive a lower bound, we write $\mathcal{T}(d_k,d_l,d_m;x,h)\ge S_1(x;h)-S_3(x;h)$ and following steps as in \cite[Lemma 3.1]{MMS}, we get $S_3(x;h) \ll_{h} (\log x)^{m-2}$.
\subsection{Estimating $S_1(x;h)$}
We need the following theorem, which we derive as a corollary of the Tauberian theorems of de la Br\'eteche \cite{RB}.
	\begin{thm}\label{three-div-function}
    As $x \to \infty$, we have 
		\begin{align*}
			\sum_{u,v,w\le x}\frac{g(u,v,w)}{[u,v,w]}d_{k-1}(u)d_{l-1}(v)d_{m-1}(w)= \nabla_{k,l,m}(h)\frac{(\log x)^{k+l+m-3}}{(k-1)!(l-1)!(m-1)!}
            +O((\log x)^{k+l+m-4}),
		\end{align*}
		where $g(u,v,w) = 1$ if the system $n\equiv-h\bmod u,n\equiv0\bmod v,n\equiv h\bmod w$ has a solution, else it is 0 and $[u,v,w]$ denotes the least common multiple of $u,v$ and $w$ and $\nabla_{k,l,m}(h)$ is the constant as in \eqref{our-conjecture-statement}.
	\end{thm}

We first complete the estimate for $S_1(x;h)$ assuming Theorem \ref{three-div-function}. Using the Chinese remainder theorem, we write
\begin{align*}
		S_1(x;h) &= \sum\limits_{{u,v,w\le x}}d_{k-1}(u)d_{l-1}(v)d_{m-1}(w) \ \sumdash_{n\le x}1\\
		&\ge\sum\limits_{{u,v,w\le x^{{1}/{3}}}}d_{k-1}(u)d_{l-1}(v)d_{m-1}(w) \ \sumdash_{n\le x}1\\ 
		&= \sum_{u,v,w\le x^{{1}/{3}}}d_{k-1}(u)d_{l-1}(v)d_{m-1}(w)g(u,v,w)\left\{\frac{x}{[u,v,w]}+O(1)\right\}\\
		&= x\sum_{u,v,w\le x^{{1}/{3}}}d_{k-1}(u)d_{l-1}(v)d_{m-1}(w)\frac{g(u,v,w)}{[u,v,w]}+O\left(\sum_{u,v,w\le x^{{1}/{3}}}d_{k-1}(u)d_{l-1}(v)d_{m-1}(w)\right).
	\end{align*}
	Using \thmref{three-div-function} we get
	$$
	S_1(x;h)\ge \frac{\nabla_{k,l,m}(h)}{3^{k+l+m-3}}\frac{x(\log x)^{k+l+m-3}}{(k-1)!(l-1)!(m-1)!}+O\left(x(\log x)^{k+l+m-4}\right).
	$$
    This completes the proof of \thmref{lower bound of triple}.
    Therefore, it remains to prove \thmref{three-div-function}.
    
\subsection{Proof of \thmref{three-div-function}}
	We consider the multiple Dirichlet series 
	$$
	F(s_1,s_2,s_3):= \sum\limits_{\substack{ u,v,w\ge 1}}\frac{g(u,v,w)}{[u,v,w]}\frac{d_{k-1}(u)}{u^{s_1}}\frac{d_{l-1}(v)}{v^{s_2}}\frac{d_{m-1}(w)}{w^{s_3}}
	$$
	on the domain $\Re(s_1),\Re(s_2),\Re(s_3)>1$.  
	For a multiplicative function $f$ of several variables, we introduce a formal Dirichlet series of several variables along with an Euler product
	$$
	\sum_{n_i\ge 1} \frac{f(n_1,\hdots,n_k)}{n_1^{s_1}\cdots n_k^{s_k}} = \prod_{p}\left(\sum_{\nu_1,\hdots,\nu_k\ge 0}\frac{f(p^{\nu_1},\hdots,p^{\nu_k})}{p^{\nu_1s_1}\cdots p^{\nu_ks_k}}\right).
	$$
	As noted in \cite{MMS}, $g$ is multiplicative and hence $F$ has an Euler product which is convergent for $\Re(s_1),\Re(s_2),\Re(s_3)>1$, namely,
	$$
	F(s_1,s_2,s_3) = \prod_{p} \left(\sum_{\nu_1,\nu_2,\nu_3\ge 0} \frac{g(p^{\nu_1},p^{\nu_2},p^{\nu_3})}{[p^{\nu_1},p^{\nu_2},p^{\nu_3}]}\frac{d_{k-1}(p^{\nu_1})d_{l-1}(p^{\nu_2})d_{m-1}(p^{\nu_3})}{p^{\nu_1s_1+\nu_2s_2+\nu_3s_3}}\right).
	$$
	We observe that if at least two of the $\nu_i$'s are $\ge1$ and $g(p^{\nu_1},p^{\nu_2},p^{\nu_3})=1$, then $p\mid 2h$. 
	Therefore, we split the Euler product into two sub-products, one for $p\mid 2h$ and the other one for $p\nmid 2h$. 
	
	Let us consider the product
	$$
	\prod_{p\ \nmid \ 2h} \left(\sum_{\nu_1,\nu_2,\nu_3\ge 0}\frac{g(p^{\nu_1},p^{\nu_2},p^{\nu_3})}{[p^{\nu_1},p^{\nu_2},p^{\nu_3}]}\frac{d_{k-1}(p^{\nu_1})d_{l-1}(p^{\nu_2})d_{m-1}(p^{\nu_3})}{p^{\nu_1s_1+\nu_2s_2+\nu_3s_3}}\right).
	$$
	We note that, the non-zero contribution in the sum is coming only from the tuples $(0,0,0)$, $(\nu_1,0,0), (0,\nu_2,0), (0,0,\nu_3)$, where $\nu_i\ge 1$ in the respective cases.
	For $(\nu_1,0,0)$, the contribution is
	$$
	\sum_{\nu_1\ge 1}\frac{d_{k-1}(p^{\nu_1})}{p^{(1+s_1)\nu_1}}. 
	$$
	From \eqref{Power series expression} we conclude that for $\Re(1+s_1)>0$ i.e., $\Re(s_1)>-1$,   
	$$
	\sum_{\nu_1\ge 1}\frac{d_{k-1}(p^{\nu_1})}{p^{(1+s_1)\nu_1}} = \left(1-\frac{1}{p^{1+s_1}}\right)^{1-k}-1.
	$$
	Similarly, for $(0,\nu_2,0)$ and $(0,0,\nu_3)$, the contributions are 
	$$
	\sum_{\nu_2\ge 1}\frac{d_{l-1}(p^{\nu_2})}{p^{(1+s_2)\nu_2}} = \left(1-\frac{1}{p^{1+s_2}}\right)^{1-l}-1,
	\text{ and }
	\sum_{\nu_3\ge 1}\frac{d_{m-1}(p^{\nu_3})}{p^{(1+s_3)\nu_3}} = \left(1-\frac{1}{p^{1+s_3}}\right)^{1-m}-1,
	$$
	for $\Re(s_2)>-1$ and $\Re(s_3)>-1$, respectively.
	Hence, for a prime number $p\nmid 2h$, the corresponding Euler factor is
	\begin{align*}
		&\left\{1+\sum_{\nu_1\ge 1}\frac{d_{k-1}(p^{\nu_1})}{p^{(1+s_1)\nu_1}}+\sum_{\nu_2\ge 1}\frac{d_{l-1}(p^{\nu_2})}{p^{(1+s_2)\nu_2}}+\sum_{\nu_3\ge 1}\frac{d_{m-1}(p^{\nu_3})}{p^{(1+s_3)\nu_3}}\right\}\\
		=&\left\{\left(1-\frac{1}{p^{1+s_1}}\right)^{1-k}+\left(1-\frac{1}{p^{1+s_2}}\right)^{1-l}+\left(1-\frac{1}{p^{1+s_3}}\right)^{1-m}-2\right\}.
	\end{align*}
	Therefore,
	\begin{align*}
		F(s_1,s_2,s_3) &= C_h(s_1,s_2,s_3)\prod_{p\ \nmid \ 2h}\left\{\left(1-\frac{1}{p^{1+s_1}}\right)^{1-k}+\left(1-\frac{1}{p^{1+s_2}}\right)^{1-l}+\left(1-\frac{1}{p^{1+s_3}}\right)^{1-m}-2\right\},
	\end{align*}
	where
	$$
	C_h(s_1,s_2,s_3) = \prod_{p \ \mid \ 2h}\left(\sum_{\nu_1,\nu_2,\nu_3\ge 0}\frac{g(p^{\nu_1},p^{\nu_2},p^{\nu_3})}{[p^{\nu_1},p^{\nu_2},p^{\nu_3}]}\frac{d_{k-1}(p^{\nu_1})d_{l-1}(p^{\nu_2})d_{m-1}(p^{\nu_3})}{p^{\nu_1s_1+\nu_2s_2+\nu_3s_3}}\right).
	$$
	
	For two infinite Euler-products $F_1(s_1,s_2,s_3)$ and $F_2(s_1,s_2,s_3)$,
    defined on a domain ${\Omega}\subseteq\mathbb{C}^3$, we denote $F_1(s_1,s_2,s_3) \approx_{\Omega} F_2(s_1,s_2,s_3)$ if ${F_1(s_1,s_2,s_3)}/{F_2(s_1,s_2,s_3)}$ is convergent on ${\Omega}$. We show that 
	$$
F(s_1,s_2,s_3)\approx_{\Omega}\zeta(1+s_1)^{k-1}\zeta(1+s_2)^{l-1}\zeta(1+s_3)^{m-1},
	$$
	on the domain
	$$
	{\Omega} := \left\{(s_1,s_1,s_3)\in \mathbb{C}^3 : \Re(s_i)>-\frac{1}{2} \  \text{for all} \  i\right\}.
	$$
	For $X,Y,Z\neq 0,1$, we define $\eta$ as follows:
	\begin{align*}
		\eta(X,Y,Z):=&\left\{\left(1-X\right)^{1-k}+\left(1-Y\right)^{1-l}+\left(1-Z\right)^{1-m}-2\right\}\left(1-X\right)^{k-1}\left(1-Y\right)^{l-1}\left(1-Y\right)^{m-1}\\
		=&\left(1-X\right)^{k-1}\left(1-Y\right)^{l-1}+\left(1-Y\right)^{l-1}\left(1-Z\right)^{m-1}+\left(1-Z\right)^{m-1}\left(1-X\right)^{k-1}\\
		&-2\left(1-X\right)^{k-1}\left(1-Y\right)^{l-1}\left(1-Z\right)^{m-1}.
	\end{align*}
	 Note that the coefficient of $X$ in $\eta(X,Y,Z)$ is $0$ and so are the coefficients of $Y,Z$ in $\eta(X,Y,Z)$. Also, the constant coefficient is $1$. Therefore $\eta(X,Y,Z)-1$ is an $\mathbb{R}$-linear combination of monomials of degree $\ge2$. In other words, $\eta(X,Y,Z)$ is a polynomial in three variables of the form 
	$$
	\eta(X,Y,Z) = 1+\sum_{2\le i_1+i_2+i_3\le k+l+m-3}a_{i_1,i_2,i_3}X^{i_1}Y^{i_2}Z^{i_3},
	$$
	where $a_{i_1,i_2,i_3}$'s belong to $\R$.
	We define $\eta_p(s_1,s_2,s_3) := \eta\left({1}/{p^{1+s_1}},{1}/{p^{1+s_2}},{1}/{p^{1+s_3}}\right)$. 
	Hence
	$$
	\eta_p(s_1,s_2,s_3)-1 = \sum_{2\le i_1+i_2+i_3\le k+l+m-3}a_{i_1,i_2,i_3}\left(\frac{1}{p^{1+s_1}}\right)^{i_1}\left(\frac{1}{p^{1+s_2}}\right)^{i_2}\left(\frac{1}{p^{1+s_2}}\right)^{i_3}. 
	$$
	The sum $\sum_{p\ \nmid \ 2h} \ \sum_{i_1+i_2+i_3=2}a_{i_1,i_2,i_3}\left(\frac{1}{p^{1+s_1}}\right)^{i_1}\left(\frac{1}{p^{1+s_2}}\right)^{i_2}\left(\frac{1}{p^{1+s_3}}\right)^{i_3}$ converges if $\Re(2+s_i+s_j)>1$, for all
    $1 \le i,j \le 3$.
	Hence $\sum_{p\ \nmid \ 2h} \ \sum_{i_1+i_2+i_3=2}a_{i_1,i_2,i_3}\left(\frac{1}{p^{1+s_1}}\right)^{i_1}\left(\frac{1}{p^{1+s_2}}\right)^{i_2}\left(\frac{1}{p^{1+s_3}}\right)^{i_3}$ is absolutely convergent in the domain 
	$$
	\Omega = \left\{(s_1,s_1,s_3)\in \mathbb{C}^3 : \Re(s_i)>-\frac{1}{2} \  \text{for all} \  i\right\}.
	$$
	Moreover, for all $(s_1,s_2,s_3)$ in $\Omega$, the sum $\sum_{p \ \nmid \ 2h}\sum_{i_1+i_2+i_3>2}a_{i_1,i_2,i_3}\left(\frac{1}{p^{1+s_1}}\right)^{i_1}\left(\frac{1}{p^{1+s_1}}\right)^{i_2}\left(\frac{1}{p^{1+s_1}}\right)^{i_3}$ is absolutely convergent.
	Therefore, the product $D_h(s_1,s_2,s_3):=\prod_{p\ \nmid \ 2h}\eta_p(s_1,s_2,s_3)$ is absolutely convergent on $\Omega$. Hence,
	\begin{align*}
		&\prod_{p \ \nmid \ 2h}\left\{\left(1-\frac{1}{p^{1+s_1}}\right)^{1-k}+\left(1-\frac{1}{p^{1+s_2}}\right)^{1-l}+\left(1-\frac{1}{p^{1+s_3}}\right)^{1-m}-2\right\}\\
        &=D_h(s_1,s_2,s_3)\times
		\prod_{p \ \nmid \ 2h}\left\{\left(1-\frac{1}{p^{1+s_1}}\right)^{1-k}
		\left(1-\frac{1}{p^{1+s_2}}\right)^{1-l}\left(1-\frac{1}{p^{1+s_3}}\right)^{1-m}\right\}
	\end{align*}
	Therefore, 
	\begin{align*}
		F(s_1,s_2,s_3) =& C_h(s_1,s_2,s_3)D_h(s_1,s_2,s_3)\zeta(1+s_1)^{k-1}\zeta(1+s_2)^{l-1}\zeta(1+s_3)^{m-1}\\
		&\times\prod_{p \ \mid \ 2h}\left\{\left(1-\frac{1}{p^{1+s_1}}\right)^{k-1}
		\left(1-\frac{1}{p^{1+s_2}}\right)^{l-1}\left(1-\frac{1}{p^{1+s_3}}\right)^{m-1}\right\}.
	\end{align*}
	
	Now we show that the function $F$ satisfies the conditions of \thmref{theorem 1, breteche} and \thmref{theorem 2, breteche}, respectively. We first note that the coefficients of $F$ are non-negative. The multiple Dirichlet series $F(s_1,s_2,s_3)$ is absolutely convergent for $\Re(s_i)>0$ for all $i$. We choose the family $\mathcal{L} :=\left\{\ell^{(1)},\hdots,\ell^{(k-1)},\ell^{(k)},\hdots,\ell^{(k+l-2)},\ell^{(k+l-1)},\hdots.,\ell^{(k+l+m-3)}\right\}$ of $n=k+l+m-3$ many non-zero linear forms defined as 
    $$
    \ell^{(a)}(s_1,s_2,s_3) = s_1, \ell^{(b)}(s_1,s_2,s_3) = s_2 \text{ and }
    \ell^{(c)}(s_1,s_2,s_3) = s_3,
    $$
    for all $1 \le a \le k-1$, $k \le b \le k+l-2$ and $k+l-1 \le c \le k+l+m-3$.
		As $\zeta(1+s)$ has only a simple pole at $s=0$, the function $H(s_1,s_2,s_3)$, defined as
		\begin{align*}
			H(\textbf{\textit{s}}) &:=F(\textbf{\textit{s}}+\pmb{\alpha})\prod_{i=1}^{k+l+m-3}\ell^{(i)}(\textbf{\textit{s}}) =F(s_1,s_2,s_3)s_1^{k-1}s_2^{l-1}s_3^{m-1}
		\end{align*}
		can be extended to the domain $\Omega$. Moreover, $H$ also satisfies the necessary growth conditions for $\delta_2={1}/{2}$. This follows from the standard convexity estimates of the Riemann zeta function (see Titchmarsh \cite[Chapter 5, eq. (5.1.4)]{Titchmarsh-book}). 
		We note that here we can take $\pmb\alpha=(0,0,0)$ and $r=|J(\pmb\alpha)|=|{j \in \{1,2,3\}:\alpha_j=0}|=3$. Hence we consider the linear forms $\{\ell^{(k+l+m-2)},\ell^{(k+l+m-1)},\ell^{(k+l+m)}\}$, defined as
		\begin{align*}
			&\ell^{(k+l+m-2)}(s_1,s_2,s_3) = e_1^*(s_1,s_2,s_3)=s_1,\\
			&\ell^{(k+l+m-1)}(s_1,s_2,s_3) = e_2^*(s_1,s_2,s_3)=s_2,\\
			&\ell^{(k+l+m)}(s_1,s_2,s_3) = e_3^*(s_1,s_2,s_3)=s_3.
		\end{align*}
Therefore, using \thmref{theorem 1, breteche}, we conclude that as $x \to \infty$,
	\begin{align*}
		\sum_{u\le x^{\beta_1}}\sum_{v\le x^{\beta_2}}\sum_{w\le x^{\beta_3}}\frac{g(u,v,w)}{[u,v,w]}d_{k-1}(u)d_{l-1}(v)d_{m-1}(w)
		&\sim x^{\langle\pmb\alpha,\pmb\beta\rangle}\left(Q_{\pmb\beta}(\log x)+O(x^{-\theta})\right),
	\end{align*}
	for some $\theta>0$, where $\pmb\beta=(\beta_1,\beta_2,\beta_3)$. The function $F$ satisfies the hypotheses of \thmref{theorem 2, breteche} for the choice $\pmb\beta=(1,1,1)$ as well. The linear form
    $$\mathcal{B}(s_1,s_2,s_3) := \left(\sum_{j=1}^{3}\beta_je_j^*\right)(s_1,s_2,s_3) = s_1+s_2+s_3$$ can be written as a $\mathbb{R}_{\ge0}$-linear combination of $\{\ell^{(i)}\}_{i=1,..,k+l+m}.$
	In fact, there exists no strict subfamily $\mathcal{L'}$ of $\{\ell^{(i)}\}_{i=1,..,k+l+m}$ such that $\mathcal{B}$ can be written as a $\mathbb{R}_{\ge0}$-linear combination of $\mathcal{L'}$ and 
		$$
		\text{card}(\mathcal{L'})-\text{rank}(\mathcal{L'}) = \text{card}\left(\{\ell^{(i)}\}_{i=1,..,k+l+m}\right) - \text{rank}\left(\{\ell^{(i)}\}_{i=1,..,k+l+m}\right).
		$$
    Hence the polynomial $Q_{\pmb\beta}$ satisfies the relation
	$$
	Q_{\pmb\beta}(\log x) = C_0x^{-\langle\pmb\alpha,\pmb\beta\rangle}I(x^{\pmb\beta})+O((\log x)^{\rho-1}),
	$$
	where $C_0:= H(0,0,0), \rho:=n+r-\text{rank}(\{\ell^{(i)}\}_{i=1}^{k+l+m}) = k+l+m-3$ 
	and 
	$$I(x^{\pmb\beta}) := \int\cdots\int_{A(x^{\pmb\beta})}\frac{dy_1\cdots dy_{k+l+m-3}}{\prod_{i=1}^{k+l+m-3} y_i^{1-\ell^{(i)}(\pmb\alpha)}},$$ 
	with 
	\begin{align*}
		A(x^{\pmb\beta}) :=&\left\{\pmb{y} \in [1,+\infty)^{k+l+m-3}: \prod_{i=1}^{k+l+m-3}y_i^{\ell^{(i)}(e_j)} \le x^{\beta_j} \ \text{ for all }j\right\}.
	\end{align*}
	Therefore,
	\begin{align*}
		I(x^{\pmb\beta}) =&\left(\int_{y_1=1}^{x}\int_{y_2=1}^{\frac{x}{y_1}}\cdots\int_{y_{k-1}=1}^{\frac{x}{y_1\cdots y_{k-2}}}\frac{dy_1\cdots dy_{k-1}}{y_1\cdots y_{k-1}}\right)\left(\int_{y_{k}=1}^{x}\int_{y_{k+1}=1}^{\frac{x}{y_k}}\cdots\int_{y_{k+l-2}=1}^{\frac{x}{y_1\cdots y_{k+l-3}}}\frac{dy_k\cdots dy_{k+l-2}}{y_k\cdots y_{k+l-2}}\right)\\
		&\times\left(\int_{y_{k+l-1}=1}^{x}\int_{y_{k+l}=1}^{\frac{x}{y_{k+l-1}}}\cdots\int_{y_{k+l+m-3}=1}^{\frac{x}{y_{k+l-1}\cdots y_{k+l+m-4}}}\frac{dy_{k+l-1}\cdots dy_{k+l+m-3}}{y_{k+l-1}\cdots y_{k+l+m-3}}\right).
    \end{align*}
    So
\begin{align*}
		I(x^{\pmb\beta}) =\frac{(\log x)^{k-1}}{(k-1)!}
        \frac{(\log x)^{l-1}}{(l-1)!}
        \frac{(\log x)^{m-1}}{(m-1)!}
		=\frac{(\log x)^{k+l+m-3}}{(k-1)!(l-1)!(m-1)!}.
	\end{align*}
	Therefore,
	\begin{align*}
		&\sum_{\substack{ u,v,w\le x}}d_{k-1}(u)d_{l-1}(v)d_{m-1}(w)\frac{g(u,v,w)}{[u,v,w]}= \nabla_{k,l,m}(h)\frac{(\log x)^{k+l+m-3}}{(k-1)!(l-1)!(m-1)!}+O((\log x)^{k+l+m-4}),
	\end{align*}
	where
	\begin{align*}
		&\nabla_{k,l,m}(h)=H(0,0,0)= C_h(0,0,0)D_h(0,0,0)\prod_{p\mid 2h}\left(1-\frac{1}{p}\right)^{k+l+m-3}.
	\end{align*} 
	Recall that,
	\begin{align*}
		D_h(0,0,0)=&\prod_{p\ \nmid \ 2h}\left\{\left(1-\frac{1}{p}\right)^{1-k}+\left(1-\frac{1}{p}\right)^{1-l}+\left(1-\frac{1}{p}\right)^{1-m}-2\right\}\left(1-\frac{1}{p}\right)^{k+l+m-3},
	\end{align*}
	and
	\begin{align*}
		C_h(0,0,0) = \prod_{p \ \mid \ 2h}\left(\sum_{\nu_1,\nu_2,\nu_3\ge 0}d_{k-1}(p^{\nu_1})d_{l-1}(p^{\nu_2})d_{m-1}(p^{\nu_3})\frac{g(p^{\nu_1},p^{\nu_2},p^{\nu_3})}{[p^{\nu_1},p^{\nu_2},p^{\nu_3}]}\right).
	\end{align*}
	Hence
	\begin{align*}
	H(0,0,0)=&\prod_{p} \left(1-\frac{1}{p}\right)^{k+l+m-3}\prod_{p\ \nmid \ 2h}\left\{\left(1-\frac{1}{p}\right)^{1-k}+\left(1-\frac{1}{p}\right)^{1-l}+\left(1-\frac{1}{p}\right)^{1-m}-2\right\}\\
	&\times\prod_{p \ \mid \ 2h}\left(\sum_{\nu_1,\nu_2,\nu_3\ge 0}d_{k-1}(p^{\nu_1})d_{l-1}(p^{\nu_2})d_{m-1}(p^{\nu_3})\frac{g(p^{\nu_1},p^{\nu_2},p^{\nu_3})}{[p^{\nu_1},p^{\nu_2},p^{\nu_3}]}\right)\\
    =&C_{k,l,m}(h)D_{k,l,m}(h)\prod_{p}\left(1-\frac{1}{p}\right)^{k+l+m-3}.
	\end{align*}
    This completes the proof of \thmref{three-div-function}.

\section{Proof of \thmref{lower bound of triple sum over primes}}
	We write the shifted convolution sum over primes as 
    $$\mathcal{T}'(d_k,d_m;x,h) =\sum_{n\le x}d_k(n+h)d_m(n-h)a(n),$$
    where
	$$
    a(n)=\begin{cases}
		1 & \text{if $n$ is prime,}\\
		0 & \text{otherwise.}\\
	\end{cases}
	$$
	Therefore,
	\begin{align*}
    \mathcal{T}'(d_k,d_m;x,h) &=\sum_{n\le x}a(n)\left(\sum_{u\mid n+h}d_{k-1}(u)\right)\left(\sum_{w\mid n-h}d_{m-1}(w)\right)\\
	&=\sum\limits_{\substack{u\le x+h\\ w\le x-h}}d_{k-1}(u)d_{m-1}(w)\sumdash_{n\le x}a(n),
	\end{align*}
	where 
	$$\sumdash_{n\le x}a(n) := \sum_{\substack{n\le x \\ n+h\equiv0 \bmod u\\ n-h\equiv0 \bmod w}}a(n) = g_1(u,w)\sum_{\substack{n\le x\\n\equiv \alpha \bmod [u,w]}}a(n),$$
	for some $0\le\alpha=\alpha(h,u,w)<[u,w]$ and $g_1(u,w) = 1$ only if the system $n\equiv-h\bmod u,n\equiv h\bmod w$ has a solution, else it is $0$. First we observe that
    $$
    \mathcal{T}'(d_k,d_m;x,h)=S_1'(x;h)+S_2'(x;h)-S_3'(x;h),
    $$
	where $S_i(x;h)$ are defined as follows:
	\begin{align*}
		&S_1'(x;h) := \sum\limits_{{u,w\le x}}d_{k-1}(u)d_{m-1}(w) \ \sumdash_{n\le x}a(n),\\
		&S_2'(x;h) := \sum\limits_{{x<u\le x+h \atop w\le x}}d_{k-1}(u)d_{m-1}(w)\ \sumdash_{n\le x}a(n),\\
		&S_3'(x;h) :=  \sum\limits_{{u\le x+h \atop x-h<w\le x}}d_{k-1}(u)d_{m-1}(w)\ \sumdash_{n\le x}a(n).
	\end{align*}
	As before, for a lower bound, we write $\mathcal{T}'(d_k,d_m;x,h)\ge S_1'(x;h)-S_3'(x;h)$.
	We also note that
    $$
    S_3'(x;h) :=  \sum\limits_{{u\le x+h\atop x-h<w\le x}}d_{k-1}(u)d_{m-1}(w)\ \sumdash_{n\le x}a(n)\le \sum\limits_{{u\le x+h\atop x-h<w\le x}}d_{k-1}(u)d_{m-1}(w)\ \sumdash_{n\le x}1,
    $$
    which is $\ll hd_k(h)(\log x)^{m-2}$ (see Lemma 3.1 in \cite{MMS} for a similar estimate). 

    \subsection{Estimating $S_1'(x;h)$} For a lower bound on $S_1'(x;h)$ we consider the following multiple Dirichlet series
    $$
    F_1(s_1,s_2):=\sum_{u,w\ge 1}\frac{g_1(u,w)}{\phi([u,w])}d_{k-1}(u)d_{m-1}(w).
    $$
    Applying \thmref{theorem 1, breteche} and \thmref{theorem 2, breteche} we prove the following theorem.
    
	\begin{thm}\label{three-div-function-prime}
    As $x \to \infty$, we have 
	$$
	\sum_{u,w\le x}\frac{g_1(u,w)}{\phi([u,w])}d_{k-1}(u)d_{m-1}(w)= \nabla_{k,m}'(h)\frac{(\log x)^{k+m-2}}{(k-1)!(m-1)!}+O((\log x)^{k+m-3}),
	$$
	where $g_1(u,w) = 1$ if the system $n\equiv-h\bmod u,n\equiv h\bmod w$ has a solution, else it is $0$. Here $[u,w]$ denotes the least common multiple of $u,w$ and $\nabla_{k,m}'(h)$ is the constant as in \eqref{constant}.
	\end{thm}
	
	Assuming Theorem \ref{three-div-function-prime}, we first complete the estimate for $S_1'(x;h)$. Let $y$ be a parameter depending on $x$ such that $1\le y\le x$ (to be chosen later). We write 
	$$
		S_1'(x;h)\ge\sum_{u,w\le y}d_{k-1}(u)d_{m-1}(w)g_1(u,w)\sum_{\substack{n\le x\\n\equiv\alpha\bmod{[u,w]}}}a(n).
    $$    
    Further restricting to those $\alpha=\alpha_h(u,w)$ that are coprime to $[u,w]$, we write
    $$
    S_1'(x;h)\ge {\mathop{{\sum}_1}}+{\mathop{{\sum}_2}}+{\mathop{{\sum}_3}},
    $$
    where
    \begin{align*}
    {\mathop{{\sum}_1}} & =\sum_{\substack{u,w\le y\\ [u,w]\le \frac{\sqrt{x}}{(\log x)^B}}}d_{k-1}(u)d_{m-1}(w)g_1(u,w)\left\{\pi\left(x;[u,w],\alpha\right)-\frac{\text{li}(x)}{\phi([u,w])}\right\},\\
	{\mathop{{\sum}_2}}& =\text{li}(x)\sum_{\substack{u,w\le y\\ [u,w]\le \frac{\sqrt{x}}{(\log x)^B}}}d_{k-1}(u)d_{m-1}(w)\frac{g_1(u,w)}{\phi([u,w])},\\
	{\mathop{{\sum}_3}}& =\sum_{\substack{u,w\le y\\ [u,w]> \frac{\sqrt{x}}{(\log x)^B}}} d_{k-1}(u)d_{m-1}(w)g_1(u,w)\pi\left(x;[u,w],\alpha\right),
	\end{align*}
	for some $B>0$, to be specified later.
   We apply \thmref{BV-a} to treat ${\mathop{{\sum}_1}}$. We write
	\begin{align*}
		\left|{\mathop{{\sum}_1}}\right| &\le \sum_{q\le \frac{\sqrt{x}}{(\log x)^B}} \max_{\gcd(a,q)=1}\left|\pi\left(x;q,a\right)-\frac{\text{li}(x)}{\phi(q)}\right|\sum_{\substack{u,w\le y\\ [u,w]=q}}d_{k-1}(u)d_{m-1}(w)g_1(u,w)\\
		&\le \sum_{q\le \frac{\sqrt{x}}{(\log x)^B}} \max_{\gcd(a,q)=1}\left|\pi\left(x;q,a\right)-\frac{\text{li}(x)}{\phi(q)}\right|\sum_{\substack{u|q\\w|q}}d_{k-1}(u)d_{m-1}(w)\\
		&\le \sum_{q\le \frac{\sqrt{x}}{(\log x)^B}}d_k(q)d_m(q) \max_{\gcd(a,q)=1}\left|\pi\left(x;q,a\right)-\frac{\text{li}(x)}{\phi(q)}\right|.
	\end{align*}
    We now choose $B>0$ such that 
    $$
    {\mathop{{\sum}_1}}\ll \frac{x}{(\log x)^A},
    $$
	for some fixed $A>0$. Next we choose $y=\frac{x^{1/4}}{(\log x)^{B/2}}$. For this choice of $y$, we also note that 
    ${\mathop{{\sum}_3}} = 0$, being an empty sum.
	For ${\mathop{{\sum}_2}}$, we apply \thmref{three-div-function-prime}. Note that
	$$
	\sum_{\substack{u,w\le \frac{x^{1/4}}{(\log x)^{B/2}}\\ [u,w]\le \frac{\sqrt{x}}{(\log x)^B}}}d_{k-1}(u)d_{m-1}(w)\frac{g_1(u,w)}{\phi([u,w])} = \sum_{\substack{u,w\le \frac{x^{1/4}}{(\log x)^{B/2}}}}d_{k-1}(u)d_{m-1}(w)\frac{g_1(u,w)}{\phi([u,w])}.
	$$
	Now, applying \thmref{three-div-function-prime}, we get 
	$$
	\text{li}(x)\sum_{\substack{u,w\le \frac{x^{1/4}}{(\log x)^{B/2}}}}d_{k-1}(u)d_{m-1}(w)\frac{g_1(u,w)}{\phi([u,w])} =	\frac{\nabla_{k,m}'(h)}{4^{k+m-2}}\frac{x(\log x)^{k+m-3}}{(k-1)!(m-1)!}
	+O(x(\log x)^{k+l+m-4}).
	$$
	Therefore it only remains to prove \thmref{three-div-function-prime} to establish \thmref{lower bound of triple sum over primes}. 
	
	\subsection{Proof of \thmref{three-div-function-prime}}
	We consider the following multiple Dirichlet series with non-negative coefficients:
	$$
	F_1(s_1,s_2) = \sum\limits_{\substack{u,w\ge 1}}\frac{g_1(u,w)}{\phi([u,w])}\frac{d_{k-1}(u)}{u^{s_1}}\frac{d_{m-1}(w)}{w^{s_2}},
	$$
	defined for $\Re(s_1),\Re(s_2)>1$.  
	As the function $g_1$ is multiplicative, we write $F$ as an Euler product, convergent for $\Re(s_1),\Re(s_2)>1$, namely,
	$$
	F_1(s_1,s_2) = \prod_{p} \left(\sum_{\nu_1,\nu_2\ge 0} \frac{g_1(p^{\nu_1},p^{\nu_2})}{\phi([p^{\nu_1},p^{\nu_2}])}\frac{d_{k-1}(p^{\nu_1})d_{m-1}(p^{\nu_2})}{p^{\nu_1s_1+\nu_2s_2}}\right).
	$$
	As before, if both $\nu_1,\nu_2 \ge1$ and $g_1(p^{\nu_1},p^{\nu_2})=1$, then $p \mid 2h$. Therefore, we split the Euler product into two sub-products, one for $p\mid 2h$ and one for $p\nmid 2h$. 
	
	First we consider the product
	$$
	\prod_{p\ \nmid \ 2h} \left(\sum_{\nu_1,\nu_2\ge 0} \frac{g_1(p^{\nu_1},p^{\nu_2})}{\phi([p^{\nu_1},p^{\nu_2}])}\frac{d_{k-1}(p^{\nu_1})d_{m-1}(p^{\nu_2})}{p^{\nu_1s_1+\nu_2s_2}}\right).
	$$
	We note that, the non-zero contributions in the sum inside the product come only from the pairs $(0,0)$, $(\nu_1,0), (0,\nu_2)$, with $\nu_1, \nu_2\ge 1$.
	For the pairs $(\nu_1,0)$, the contribution is
	$$
	\left(1-\frac{1}{p}\right)^{-1}\sum_{\nu_1\ge 1}\frac{d_{k-1}(p^{\nu_1})}{p^{(1+s_1)\nu_1}}. 
	$$
    Clearly the infinite series converges for $\Re(1+s_1)>0$, i.e., $\Re(s_1)>-1$ and by \eqref{Power series expression} we have
	$$
	\left(1-\frac{1}{p}\right)^{-1}\sum_{\nu_1\ge 1}\frac{d_{k-1}(p^{\nu_1})}{p^{(1+s_1)\nu_1}} =	\left(1-\frac{1}{p}\right)^{-1}\left\{ \left(1-\frac{1}{p^{1+s_1}}\right)^{1-k}-1\right\}.
	$$
	Similarly, for $\Re(s_2)>-1$, the contribution from the pairs $(0,\nu_2)$ is
	$$
	\left(1-\frac{1}{p}\right)^{-1}\sum_{\nu_2\ge 1}\frac{d_{m-1}(p^{\nu_2})}{p^{(1+s_2)\nu_2}} =
	\left(1-\frac{1}{p}\right)^{-1}\left\{ \left(1-\frac{1}{p^{1+s_2}}\right)^{1-m}-1\right\}.
	$$
	Hence, for $p\nmid 2h$, the corresponding Euler factor is 
	\begin{align*}
		&1+\left(1-\frac{1}{p}\right)^{-1}\left\{\left(1-\frac{1}{p^{1+s_1}}\right)^{1-k}+\left(1-\frac{1}{p^{1+s_2}}\right)^{1-m}-2\right\}.
	\end{align*}
	Therefore,
	\begin{align*}
		F_1(s_1,s_2) = C_h'(s_1,s_2)\prod_{p \ \nmid \ 2h}\left( 1+\left(1-\frac{1}{p}\right)^{-1}\left\{\left(1-\frac{1}{p^{1+s_1}}\right)^{1-k}+\left(1-\frac{1}{p^{1+s_2}}\right)^{1-m}-2\right\}\right),
	\end{align*}
	where
	$$
	C_h'(s_1,s_2) = \prod_{p\ \mid \ 2h}\left(\sum_{\nu_1,\nu_2\ge 0}\frac{g_1(p^{\nu_1},p^{\nu_2})}{\phi([p^{\nu_1},p^{\nu_2}])}\frac{d_{k-1}(p^{\nu_1})d_{m-1}(p^{\nu_2})}{p^{\nu_1s_1+\nu_2s_2}}\right).
	$$
    We show that $F_1(s_1,s_2)$ is same as $\zeta(1+s_1)^{k-1}\zeta(1+s_2)^{m-1}$, up to an infinite product which is convergent on the domain ${\Omega}\subseteq\mathbb{C}^2$ defined by
	$$
	{\Omega} := \left\{(s_1,s_2)\in \mathbb{C}^2 : \Re(s_i)>-\frac{1}{2} \  \text{for all} \  i\right\}.
	$$
	For $X,Y\neq 0,1$, we define $\eta'$ as follows:
	\begin{align*}
		\eta'(X,Y):&=\left[1+c\left\{\left(1-X\right)^{1-k}+\left(1-Y\right)^{1-m}-2\right\}\right]\left(1-X\right)^{k-1}\left(1-Y\right)^{m-1}\\
		&=c\left(1-Y\right)^{m-1}+c\left(1-X\right)^{k-1}
		+(1-2c)\left(1-X\right)^{k-1}\left(1-Y\right)^{m-1}.
	\end{align*}
	The coefficient of $X$ in $\eta'(X,Y)$ is $-c(k-1)-(1-2c)(k-1) = (c-1)(k-1)$. Similarly, the coefficient of $Y$ in $\eta'(X,Y)$ is $(c-1)(m-1)$. So, $\eta'(X,Y)$ is a polynomial of the form 
	\begin{align*}
		\eta'(X,Y) = &1+(c-1)(k-1)X+(c-1)(m-1)Y+\sum_{i_1, i_2 \ge 0 \atop 2\le i_1+i_2\le k+m-2}a_{i_1,i_2}X^{i_1}Y^{i_2},
	\end{align*}
	for some $a_{i_1,i_2} \in \R$.

    For every prime $p \nmid 2h$, let $\eta_p(s_1,s_2):=\eta'\left(\frac{1}{p^{1+s_1}},\frac{1}{p^{1+s_2}}\right)$, together with $c=p/(p-1)$. So we have
    \begin{align*}
    \eta_p(s_1,s_2)-1= \frac{(k-1)}{(p-1)p^{1+s_1}}+\frac{(m-1)}{(p-1)p^{1+s_2}}+
		\sum_{i_1, i_2 \ge 0 \atop 2\le i_1+i_2\le k+m-2}a_{i_1,i_2}\left(\frac{1}{p^{1+s_1}}\right)^{i_1}\left(\frac{1}{p^{1+s_2}}\right)^{i_2}.
    \end{align*}
    Note that the sum $\sum_{p \ \nmid \ 2h}\frac{1}{(p-1)p^{1+s_1}}$ converges for $\Re(2+s_1)>1$, i.e., for $\Re(s_1)>-1$. Also, the sum
    $$\sum_{p\ \nmid \ 2h} \ \sum_{i_1, i_2 \ge 0 \atop 2\le i_1+i_2\le k+m-2}a_{i_1,i_2}\left(\frac{1}{p^{1+s_1}}\right)^{i_1}\left(\frac{1}{p^{1+s_2}}\right)^{i_2}$$ converges if
    $$\Re(2+2s_1)>1, \Re(2+2s_2)>1 \ \text{ and } \ \Re(2+s_1+s_2)>1.$$
	Hence all these sums are convergent in the domain 
	$$
	\Omega = \left\{(s_1,s_2)\in \mathbb{C}^2 : \Re(s_i)>-\frac{1}{2} \  \text{for all} \  i\right\}.
	$$
	Therefore, the product $D_h'(s_1,s_2):=\prod_{p\ \nmid \ 2h}\eta_p'(s_1,s_2)$ is absolutely convergent on $\Omega$. Hence,
	\begin{align*}
		&\prod_{p\ \nmid \ 2h}\left( 1+\left(1-\frac{1}{p}\right)^{-1}\left\{\left(1-\frac{1}{p^{1+s_1}}\right)^{1-k}+\left(1-\frac{1}{p^{1+s_2}}\right)^{1-m}-2\right\}\right)\\
		=&D_h'(s_1,s_2)\times
		\prod_{p\ \nmid \ 2h}\left(\left(1-\frac{1}{p^{1+s_1}}\right)^{1-k}
		\left(1-\frac{1}{p^{1+s_2}}\right)^{1-m}\right)
	\end{align*}
	Therefore, 
	\begin{align*}
		F_1(s_1,s_2) = & \ C_h'(s_1,s_2)D_h'(s_1,s_2)\zeta(1+s_1)^{k-1}\zeta(1+s_2)^{m-1}\\
		&\times \prod_{p\ \mid \ 2h}\left(\left(1-\frac{1}{p^{1+s_1}}\right)^{k-1}
		\left(1-\frac{1}{p^{1+s_2}}\right)^{m-1}\right).
	\end{align*}
	
    Now we show that the function $F_1$ satisfies the conditions of \thmref{theorem 1, breteche} and \thmref{theorem 2, breteche} respectively. From the above expression, we see that $F_1(s_1,s_2)$ is absolutely convergent if $\Re(s_1)>0$ and $\Re(s_2)>0$. So we take $\pmb\alpha=(0,0)$. We choose the family $\mathcal{L} :=\left\{\ell^{(1)},\hdots,\ell^{(k-1)},\ell^{(k)},\hdots,\ell^{(k+m-2)}\right\}$ of $n=k+m-2$ many non-zero linear forms defined as 
    $$
    \ell^{(a)}(s_1,s_2) = s_1, \text{ and }
    \ell^{(b)}(s_1,s_2) = s_2,
    $$
    for all $1 \le a \le k-1$ and $k \le b \le k+m-2$.
	As $\zeta(1+s)$ has only a simple pole at $s=0$, the function $H_1(s_1,s_2)$, defined as
	\begin{align*}
		H_1(\textbf{\textit{s}}) &:=F_1(\textbf{\textit{s}}+\pmb{\alpha})\prod_{i=1}^{k+m-2}\ell^{(i)}(\textbf{\textit{s}}) =F_1(s_1,s_2)s_1^{k-1}s_2^{m-1}
	\end{align*}
	can be extended to the domain $\Omega$. Therefore, we choose $\delta_1={1}/{2}$. Moreover, $F_1$ also satisfies the necessary growth conditions for $\delta_2={1}/{2}$. This follows from the standard convexity estimates of the Riemann zeta function (see Titchmarsh \cite[Chapter 5, eq. (5.1.4)]{Titchmarsh-book}). 
	
	
	As $\pmb\alpha=(0,0)$, we have $r=|J(\pmb\alpha)|=|\{j \in \{1,2\}:\alpha_j=0\}|=2$. Hence, we consider the linear forms $\{\ell^{(k+m-1)},\ell^{(k+m)}\}$ defined as
	\begin{align*}
		&\ell^{(k+m-1)}(s_1,s_2) = e_1^*(s_1,s_2)=s_1,\\
		&\ell^{(k+m)}(s_1,s_2) = e_2^*(s_1,s_2)=s_2.
	\end{align*}
	Therefore, using \thmref{theorem 1, breteche}, we conclude that as $x \to \infty$,
	\begin{align*}
		\sum_{u\le x^{\beta_1}}\sum_{w\le x^{\beta_2}}\frac{g_1(u,w)}{\phi([u,w])}d_{k-1}(u)d_{m-1}(w)
		&\sim x^{\langle\pmb\alpha,\pmb\beta\rangle}\left(Q_{\pmb\beta}(\log x)+O(x^{-\theta})\right),
	\end{align*}
	for some $\theta>0.$ The function $F_1$ satisfies the hypotheses of \thmref{theorem 2, breteche} for the choice $\pmb\beta=(1,1)$ as well. The linear form $\mathcal{B}(s_1,s_2) := \left(\sum_{j=1}^{2}\beta_je_j^*\right)(s_1,s_2) = s_1+s_2$ can be written as an $\mathbb{R}_{\ge0}$-linear combination of $\{\ell^{(i)}\}_{i=1,..,k+m}.$
	In fact, there exists no strict subfamily $\mathcal L'$ of $\{\ell^{(i)}\}_{i=1,..,k+m}$ such that $\mathcal{B}$ can be written as an $\mathbb{R}_{\ge0}$-linear combination of $\mathcal{L'}$ and 
	$$
	\text{card}(\mathcal{L'})-\text{rank}(\mathcal{L'}) = \text{card}\left(\{\ell^{(i)}\}_{i=1,..,k+m}\right) - \text{rank}\left(\{\ell^{(i)}\}_{i=1,..,k+m}\right).
	$$
	Hence the polynomial $Q_{\pmb\beta}$ satisfies the relation
	$$
	Q_{\pmb\beta}(\log x) = H_1(0,0)x^{-\langle\pmb\alpha,\pmb\beta\rangle}I(x^{\pmb\beta})+O((\log x)^{\rho-1}),
	$$
	where $\rho:=n+r-\text{rank}(\{\ell^{(i)}\}_{i=1}^{k+m}) = k+m-2$ 
	and 
	$$I(x^{\pmb\beta}) := \int\cdots\int_{A(x^{\pmb\beta})}\frac{dy_1\cdots dy_{k+m-2}}{\prod_{i=1}^{k+m-2} y_i^{1-\ell^{(i)}(\pmb\alpha)}},$$ 
	with 
	\begin{align*}
		A(x^{\pmb\beta}) :=&\left\{\pmb{y} \in [1,+\infty)^{k+m-2}: \prod_{i=1}^{k+m-2}y_i^{\ell^{(i)}(e_j)} \le x^{\beta_j} \ \ \text{ for all } j\right\}.
	\end{align*}
	Therefore, as before
	\begin{align*}
		I(x^{\pmb\beta}) =&\left(\int_{y_1=1}^{x}\int_{y_2=1}^{\frac{x}{y_1}}\cdots\int_{y_{k-1}=1}^{\frac{x}{y_1\cdots y_{k-2}}}\frac{dy_1\cdots dy_{k-1}}{y_1\cdots y_{k-1}}\right)\left(\int_{y_{k}=1}^{x}\int_{y_{k+1}=1}^{\frac{x}{y_k}}\cdots\int_{y_{k+m-2}=1}^{\frac{x}{y_1\cdots y_{k+m-3}}}\frac{dy_k\cdots dy_{k+m-2}}{y_k\cdots y_{k+m-2}}\right)\\
		=&\frac{(\log x)^{k-1}}{(k-1)!}\frac{(\log x)^{m-1}}{(m-1)!}
		=\frac{(\log x)^{k+m-2}}{(k-1)!(m-1)!}.
	\end{align*}
	Therefore,
	\begin{align*}
		&\sum_{\substack{u,w\le x}}d_{k-1}(u)d_{m-1}(w)\frac{g_1(u,w)}{\phi([u,w])}= \nabla_{k,m}'(h)\frac{(\log x)^{k+m-2}}{(k-1)!(m-1)!}+O((\log x)^{k+m-3}),
	\end{align*}
	where
	\begin{align*}
		&\nabla_{k,m}'(h)=H_1(0,0)= C_h'(0,0)D_h'(0,0)\prod_{p \ \mid \ 2h}\left(1-\frac{1}{p}\right)^{k+m-2}.
	\end{align*} 
	We recall that
	\begin{align*}
		D_h'(0,0)=&\prod_{p\ \nmid \ 2h}\left(1+\left(1-\frac{1}{p}\right)^{-1}\left\{\left(1-\frac{1}{p}\right)^{1-k}+\left(1-\frac{1}{p}\right)^{1-m}-2\right\}\right)\left(1-\frac{1}{p}\right)^{k+m-2}\\
		=&\prod_{p\ \nmid \ 2h}\left(1-2\left(1-\frac{1}{p}\right)^{-1}+\left(1-\frac{1}{p}\right)^{-k}+\left(1-\frac{1}{p}\right)^{-m}\right)\left(1-\frac{1}{p}\right)^{k+m-2}\\
        =&D_{k,m}'(h)\prod_{p\ \nmid \ 2h}\left(1-\frac{1}{p}\right)^{k+m-2}
	\end{align*}
	and
	\begin{align*}
		C_h'(0,0) = \prod_{p \ \mid \ 2h}\left(\sum_{\nu_1,\nu_2\ge 0}d_{k-1}(p^{\nu_1})d_{m-1}(p^{\nu_2})\frac{g_1(p^{\nu_1},p^{\nu_2})}{\phi([p^{\nu_1},p^{\nu_2}])}\right).
	\end{align*}
	This completes the proof of \thmref{lower bound of triple sum over primes}.

	\subsection{Proof of Remark \ref{rmk}}
    From \lemref{CRT}, we note that $g_1(p^{\nu_1},p^{\nu_2})=1$ if $\gcd(p^{\nu_1},p^{\nu_2})|2h$, otherwise $g_1(p^{\nu_1},p^{\nu_2})=0$. Let $p$ be an odd prime number and $v_p(h)=\alpha$. Then at most one of $\nu_1,\nu_2$ can be $>\alpha$. If $\nu_1 >\alpha$, the contribution is
    \begin{align*}
    \sum_{\nu_1>\alpha}\sum_{\nu_2=0}^{\alpha}\frac{d_{k-1}(p^{\nu_1})d_{m-1}(p^{\nu_2})}{\phi(p^{\nu_1})}&=\left(1-\frac{1}{p}\right)^{-1}\sum_{\nu_1>\alpha}\sum_{\nu_2=0}^{\alpha}\frac{d_{k-1}(p^{\nu_1})d_{m-1}(p^{\nu_2})}{p^{\nu_1}}\\
    &=\left(1-\frac{1}{p}\right)^{-1}d_{m}(p^{\alpha})\sum_{\nu_1>\alpha}\frac{d_{k-1}(p^{\nu_1})}{p^{\nu_1}}.
    \end{align*}
    Using \eqref{h}, we see that the contribution for $\nu_1>\alpha$ is
$$
d_{m}(p^{\alpha})\sum_{\nu_1>\alpha}\frac{d_{k}(p^{\nu_1})}{p^{\nu_1}}-\left(1-\frac{1}{p}\right)^{-1}\frac{d_{k}(p^{\alpha})d_m(p^{\alpha})}{p^{\alpha+1}}.
$$
Similarly, if $\nu_2>\alpha$, the contribution is 
$$
d_{k}(p^{\alpha})\sum_{\nu_2>\alpha}\frac{d_{m}(p^{\nu_2})}{p^{\nu_2}}-\left(1-\frac{1}{p}\right)^{-1}\frac{d_{k}(p^{\alpha})d_{m}(p^{\alpha})}{p^{\alpha+1}}.
$$
Now we consider the case when $\nu_1,\nu_2\le\alpha$. For $0\le\nu_1=\nu_2\le\alpha$, we get
$$
\sum_{\nu=0}^{\alpha}\frac{d_{k-1}(p^{\nu})d_{m-1}(p^{\nu})}{\phi(p^{\nu})}=1+\left(1-\frac{1}{p}\right)^{-1}\sum_{\nu=1}^{\alpha}\frac{d_{k-1}(p^{\nu})d_{m-1}(p^{\nu})}{p^{\nu}}.
$$
For $0\le\nu_2<\nu_1\le\alpha$, we get
$$
\sum_{\nu_1=1}^{\alpha}\sum_{\nu_2=0}^{\nu_1-1}\frac{d_{k-1}(p^{\nu_1})d_{m-1}(p^{\nu_2})}{\phi(p^{\nu_1})}=\left(1-\frac{1}{p}\right)^{-1}\sum_{\nu_1=1}^{\alpha}\frac{d_{k-1}(p^{\nu_1})d_{m}(p^{\nu_1-1})}{p^{\nu_1}}.
$$
For $0\le\nu_1<\nu_2\le\alpha$, we get
$$
\sum_{\nu_2=1}^{\alpha}\sum_{\nu_1=0}^{\nu_2-1}\frac{d_{k-1}(p^{\nu_1})d_{m-1}(p^{\nu_2})}{\phi(p^{\nu_2})}= \left(1-\frac{1}{p}\right)^{-1}\sum_{\nu_2=1}^{\alpha}\frac{d_k(p^{\nu_2-1})d_{m-1}(p^{\nu_2})}{p^{\nu_2}}.
$$
Therefore, by \eqref{recursive-relation}, the contribution from the indices $\nu_1,\nu_2\le\alpha$ is
$$
1+\left(1-\frac{1}{p}\right)^{-1}\sum_{\nu=1}^{\alpha}\frac{d_k(p^{\nu})d_m(p^{\nu})-d_k(p^{\nu-1})d_m(p^{\nu-1})}{p^{\nu}}.
$$
Now let $\nu_2(h)=\alpha$, i.e., $\nu_2(2h)=\alpha+1$. We observe that $g_1(2^{\nu_1},2^{\nu_2})=1$ if $\gcd(2^{\nu_1},2^{\nu_2})|2h$, i.e., $\min(\nu_1,\nu_2)\le\alpha+1$. When $\nu_1>\alpha+1$, the contribution is
$$
d_{m}(2^{\alpha+1})\sum_{\nu_1>\alpha+1}\frac{d_{k}(2^{\nu_1})}{2^{\nu_1}}-\frac{d_{k}(2^{\alpha+1})d_{m}(2^{\alpha+1})}{2^{\alpha+1}}.
$$
Similarly, for $\nu_2>\alpha+1$, the contribution is 
$$
d_{k}(2^{\alpha+1})\sum_{\nu_2>\alpha+1}\frac{d_{m}(2^{\nu_2})}{2^{\nu_2}}-\frac{d_{m}(2^{\alpha+1})d_{k}(2^{\alpha+1})}{2^{\alpha+1}}.
$$
When both $\nu_1$ and $\nu_2$ are $\le\alpha+1$, the contribution is 
$$
1+\sum_{\nu=1}^{\alpha+1}\frac{d_k(2^{\nu})d_m(2^{\nu})-d_k(2^{\nu-1})d_m(2^{\nu-1})}{2^{\nu-1}}.
$$
This completes the proof of Remark \ref{rmk}.

\begin{rmk}\rm
    It is worth noting that for the case $k=m=2$, it is possible to improve the lower bound in \thmref{lower bound of triple sum over primes}. The usual divisor function $d(n)$ can be written as
    \begin{equation}\label{DH-method}
    d(n) = 2\sum_{\substack{u|n\\u\le \sqrt{n}}}1 -
    \begin{cases}
        1 & \text{  if $n$ is a square},\\
        0 & \text{ otherwise}.
    \end{cases}
    \end{equation}
    Applying \eqref{DH-method}, together with \thmref{three-div-function-prime} and Theorem \ref{BV-a}, the factor ${1}/{4^2}$ improves to ${1}/{4}$.
\end{rmk}

\section{Getting to $\nabla_{k,l,m}(h)$ using probabilistic route}
We have already seen in \thmref{three-div-function} how $\nabla_{k,l,m}(h)$ emerges from the theory of multiple Dirichlet series as an application of the Tauberian theorems due to de la Bret\`eche \cite{RB}. We now briefly indicate how to get to $\nabla_{k,l,m}(h)$ using a probabilistic route. 

Following \cite{Ng}, let $\{X_p\}_p,\{Y_p\}_p,\{Z_p\}_p$ be sequences of random variables, indexed by the prime numbers defined as follows:
	 $$
	 X_p(n) = d_k\left(p^{v_p(n+h)}\right), \quad Y_p(n) = d_l\left(p^{v_p(n)}\right), \quad Z_p(n) = d_m\left(p^{v_p(n-h)}\right),
	 $$
	 where $v_p(n)$ denotes the $p$-adic valuation of $n$.
We define three more random variables $X,Y$ and $Z$ as follows:
	 $$
	 X(n)=\prod_{p}X_p(n), \quad Y(n)=\prod_{p}Y_p(n), \quad Z(n)=\prod_{p}Z_p(n).
	 $$
We recall that for the random variable $X:\N\to\C$, its expectation is given by
	 $$
	 \E(X) = \sum_{i\in \text{Im}(X)} i \cdot \mathbb{P}(X=i)
	 $$
	 where $\text{Im}(X)=\{X(n) : n\in\N\}$ and for any subset $B\subseteq \N$, we consider the following natural definition
	 $$
	 \mathbb{P}(B) := \lim\limits_{x\to\infty} \frac{\#\left\{1\le n\le x| n\in B\right\}}{x}.
	 $$
It is easy to deduce that as $x \to \infty$,
	 $$
	 \E(Y)\sim \frac{1}{x}\sum_{n\le x}d_l(n).
	 $$
Similarly, $\E(X),\E(Z)$ give the average order of $d_k,d_m$, respectively, for $h\le x^{1-\epsilon}$ for some $0<\epsilon<1$. It can be noted that $X,Y,Z$ are not mutually independent random variables. And based on these considerations, Conjecture \ref{tao's conjecture} (see \cite{Tao},\cite{Ng}) can be written as
	 $$
	 \frac{1}{x}\sum_{n\le x}d_{k}(n+h)d_{l}(n)\sim\left\{\prod_{p \text{ prime}}\frac{\E(X_pY_p)}{\E(X_p)\E(Y_p)}\right\}\left(\frac{1}{x}\sum_{n\le x}d_{k}(n+h)\right)\left(\frac{1}{x}\sum_{n\le x}d_{l}(n)\right),
	 $$  
	 for $h\le x^{1-\epsilon}$ and $x\to\infty$. It is therefore reasonable to expect the following triple convolution analogue of the above conjecture:  
	 
	 \begin{conj}\label{conj-exp}
     Let $\epsilon>0$ and integers $k,l,m\ge 2$. Then for $0<h\le x^{1-\epsilon}$, as $x \to \infty$
	 	\begin{equation}\label{Probabilistic-approach-to-conjecture}
	 	\sum_{n\le x}d_k(n+h)d_l(n)d_m(n-h)\sim\left( \prod_{p \text{ prime}}\frac{\E(X_pY_pZ_p)}{\E(X_p)\E(Y_p)\E(Z_p)}\right)\frac{x(\log x)^{k+l+m-3}}{(k-1)!(l-1)!(m-1)!}.
	 	\end{equation}
	 \end{conj}

The verification of the fact that
$$
\nabla_{k,l,m}(h)=\prod_{p \text{ prime}}\frac{\E(X_pY_pZ_p)}{\E(X_p)\E(Y_p)\E(Z_p)},
$$
involves calculating the required expectations. To this end, it is possible to derive
the following:

\begin{lem}\label{Expectations calculation} 
Let $p$ be a prime number.
	 	\begin{enumerate}
	 		\item[(i)] Then we have
	 		$$
	 		\E(X_p)=\left(1-\frac{1}{p}\right)^{1-k};\quad \E(Y_p)=\left(1-\frac{1}{p}\right)^{1-l};\quad \E(Z_p)=\left(1-\frac{1}{p}\right)^{1-m}.
	 		$$
	 		\item[(ii)] If $p$ is an odd prime such that $p\nmid h$, then
	 		$$
	 		\E(X_pY_pZ_p) = \left(1-\frac{1}{p}\right)^{1-k}+\left(1-\frac{1}{p}\right)^{1-l}+\left(1-\frac{1}{p}\right)^{1-m}-2.
	 		$$
	 		\item[(iii)]\label{last case-Expectations calculation} If 
            $p$ is an odd prime such that $v_p(h)=\alpha>0$, then
	 		\begin{align*}
	 		\E(X_pY_pZ_p) = & \left(1-\frac{1}{p}\right)\left(\sum_{i=0}^{\alpha-1}\frac{d_k(p^i)d_l(p^i)d_m(p^i)}{p^i}+\sum_{i=1}^{\infty}\frac{d_k(p^{\alpha+i})d_l(p^i)d_m(p^i)}{p^{\alpha+i}}\right.\\
            &\left.+\sum_{i=1}^{\infty}\frac{d_k(p^{i})d_l(p^{\alpha+i})d_m(p^i)}{p^{\alpha+i}}	+\sum_{i=1}^{\infty}\frac{d_k(p^{i})d_l(p^{i})d_m(p^{\alpha+i})}{p^{\alpha+i}}\right)\\
            &+\left(1-\frac{3}{p}\right)\frac{d_k(p^{\alpha})d_l(p^{\alpha})d_m(p^{\alpha})}{p^{\alpha}}.
	 		\end{align*}
	 		\item[(iv)] If $v_2(h)=\alpha \ge 0$, then
	 		\begin{align*}
	 			\E(X_2Y_2Z_2)= &\frac{1}{2}\sum_{i=0}^{\alpha-1}\frac{d_k(2^i)d_l(2^i)d_m(2^i)}{2^i}+\frac{1}{2}\sum_{i=1}^{\infty}\frac{d_k(2^{\alpha+i})d_l(2^{\alpha})d_m(2^{\alpha+1})}{2^{\alpha+i}}\\
	 			&+\frac{1}{2}\sum_{i=1}^{\infty}\frac{d_k(2^{\alpha})d_l(2^{\alpha+i})d_m(2^{\alpha})}{2^{\alpha+i}}+\frac{1}{2}\sum_{i=1}^{\infty}\frac{d_k(2^{\alpha+1})d_l(2^{\alpha})d_m(2^{\alpha+i})}{2^{\alpha+i}}\\
	 			&-\frac{d_k\left(2^{\alpha+1}\right)d_l\left(2^{\alpha}\right)d_m\left(2^{\alpha+1}\right)}{2^{\alpha+1}}.
	 		\end{align*}
	 	\end{enumerate}	
	 \end{lem}

This lemma therefore implies the following restatement of Conjecture \ref{conj-exp}.

\begin{conj}
Let $\epsilon>0$ and integers $k,l,m\ge 2$. Then for $0<h\le x^{1-\epsilon}$, as $x \to \infty$,
	\begin{equation*}
	\sum_{n\le x}d_k(n+h)d_l(n)d_m(n-h)\sim \varrho_{k,l,m}(h)\psi_{k,l,m}(h)\prod_{p}\left(1-\frac{1}{p}\right)^{k+l+m-3}\frac{x(\log x)^{k+l+m-3}}{(k-1)!(l-1)!(m-1)!},
	\end{equation*}
	 where
	 $$
	 \varrho_{k,l,m}(h)=\prod_{p \ \nmid \ h}\left\{\left(1-\frac{1}{p}\right)^{1-k}+\left(1-\frac{1}{p}\right)^{1-l}+\left(1-\frac{1}{p}\right)^{1-m}-2\right\},
	 $$
	 and $\psi_{k,l,m}(h)$ is defined multiplicatively as
	 \begin{align*}
	 	\psi_{k,l,m}(p^{\alpha})= & \left(1-\frac{1}{p}\right)\left(\sum_{i=0}^{\alpha-1}\frac{d_k(p^i)d_l(p^i)d_m(p^i)}{p^i}+\sum_{i=1}^{\infty}\frac{d_k(p^{\alpha+i})d_l(p^i)d_m(p^i)}{p^{\alpha+i}}\right.\\
	 	&\left.+\sum_{i=1}^{\infty}\frac{d_k(p^{i})d_l(p^{\alpha+i})d_m(p^i)}{p^{\alpha+i}}	+\sum_{i=1}^{\infty}\frac{d_k(p^{i})d_l(p^{i})d_m(p^{\alpha+i})}{p^{\alpha+i}}\right)\\
        & +\left(1-\frac{3}{p}\right)\frac{d_k(p^{\alpha})d_l(p^{\alpha})d_m(p^{\alpha})}{p^{\alpha}}
	 \end{align*}
	 for odd prime numbers $p$ and 
	 \begin{align*}
	 	\psi_{k,l,m}(2^{\alpha})= &\frac{1}{2}\sum_{i=0}^{\alpha-1}\frac{d_k(2^i)d_l(2^i)d_m(2^i)}{2^i}+\frac{1}{2}\sum_{i=1}^{\infty}\frac{d_k(2^{\alpha+i})d_l(2^{\alpha})d_m(2^{\alpha+1})}{2^{\alpha+i}}\\
	 	+&\frac{1}{2}\sum_{i=1}^{\infty}\frac{d_k(2^{\alpha})d_l(2^{\alpha+i})d_m(2^{\alpha})}{2^{\alpha+i}}+\frac{1}{2}\sum_{i=1}^{\infty}\frac{d_k(2^{\alpha+1})d_l(2^{\alpha})d_m(2^{\alpha+i})}{2^{\alpha+i}}\\
	 	-&\frac{d_k\left(2^{\alpha+1}\right)d_l\left(2^{\alpha}\right)d_m\left(2^{\alpha+1}\right)}{2^{\alpha+1}}.
	 \end{align*}
	 \end{conj}

Thus, we are left with the verification of the fact that $C_{k,l,m}(h)=\psi_{k,l,m}(h)$,
where $C_{k,l,m}(h)$ is as in Conjecture \ref{conj-nabla}. This verification and proof of \lemref{Expectations calculation} take a number of pages and hence we only give some key steps.

\subsection{Proof of \lemref{Expectations calculation}}
Part (i) can already be found in \cite[Lemma 4.2]{Ng}. For part (ii), we note that
$p$ can divide at most one of $n,n\pm h$. So we will have four possible choices for
$(X_p(n),Y_p(n),Z_p(n))$, namely, $(1,1,1)$, $(d_k(p^i),1,1)$, $(1,d_l(p^i),1)$ and
$(1,1, d_m(p^i))$ for $i \ge 1$, where $i$ denotes the $p$-adic valuation of $n+h$, or $n$, or $n-h$, respectively.
Computing the required probability, we get the desired result.

Part (iii) is more involved. If we assume $p \nmid n$, then we get $(X_p(n),Y_p(n),Z_p(n))=(1,1,1)$ and hence the contribution to $\E(X_pY_pZ_p)$ is
$1-1/p$. Now we consider $v_p(n)=i$ with $i\ge 1$. 
In order to determine $v_p(n+h)$ and $v_p(n-h)$, we write $n=p^in'$ and $h=p^{\alpha}h'$ with $\gcd(n',p)=\gcd(h',p)=1$. Therefore,
$$
n\pm h = p^{\text{min}(i,\alpha)}\left(n'p^{i-\text{min}(i,\alpha)}\pm h'p^{\alpha-\text{min}(i,\alpha)}\right).
$$ 
Hence, for the case $i\neq\alpha$, the contribution to $\E(X_pY_pZ_p)$ is 
$$
\sum_{i=1\atop i\neq\alpha}^{\infty} d_k\left(p^{\text{min}(i,\alpha)}\right)d_l(p^i)d_m\left(p^{\text{min}(i,\alpha)}\right)\left(\frac{1}{p^i}-\frac{1}{p^{i+1}}\right).
$$
Now if $i=\alpha$, then $v_p(n \pm h)=\alpha+v_p(n'\pm h')$. 
We first observe that, at most one of $v_p(n \pm h)$ can be more than $\alpha$.
In particular, if $v_p(n \pm h)=\alpha$, i.e., 
$n'\not\equiv 0,h',-h' \bmod p$, we get the contribution
$$
\frac{d_k(p^{\alpha})d_l(p^{\alpha})d_m(p^{\alpha})}{p^{\alpha}}\left(1-\frac{3}{p}\right).
$$
Now for the case $v_p(n + h)>\alpha$ (or, $v_p(n - h)>\alpha$), we consider the following events for $j \ge 1$:
$$
A_j^{+} = \left\{n\in\N : v_p(n)=\alpha, \ v_p(n'+ h')=j\right\}
\ \left(\text{resp. } A_j^{-} = \left\{n\in\N : v_p(n)=\alpha, \ v_p(n'- h')=j\right\}\right).
$$ 
Then we have (see page 134 in \cite{Ng}) 
$$
\mathbb{P}(A_j^{\pm})=\frac{1}{p^{\alpha}}\left(\frac{1}{p^j}-\frac{1}{p^{j+1}}\right).
$$
Therefore, contributions for the cases $v_p(n + h)>\alpha$ and $v_p(n - h)>\alpha$ are
$$
\sum_{j=1}^{\infty}\frac{d_k\left(p^{\alpha+j}\right)d_l\left(p^{\alpha}\right)d_m\left(p^{\alpha}\right)}{p^{\alpha}}\left(\frac{1}{p^j}-\frac{1}{p^{j+1}}\right)
\text{ and }
\sum_{j=1}^{\infty}\frac{d_k\left(p^{\alpha}\right)d_l\left(p^{\alpha}\right)d_m\left(p^{\alpha+j}\right)}{p^{\alpha}}\left(\frac{1}{p^j}-\frac{1}{p^{j+1}}\right),
$$
respectively. Putting all these together, we get the desired result.

For part (iv), the computation when $h$ is odd is simple and the contribution to $\E(X_2Y_2Z_2)$ turns out to be
$$
\sum_{i=1}^{\infty}\frac{d_l(2^i)}{2^{i+1}},
$$ 
when $n$ is even and
$$
m\sum_{i=1}^{\infty}\frac{d_k(2^i)}{2^{i+1}}+k\sum_{i=1}^{\infty}\frac{d_m(2^i)}{2^{i+1}}-\frac{km}{2},
$$
when $n$ is odd. So this proves the formula for $\alpha=0$. So let $\alpha>0$. If $n$ is odd, the contribution to $\E(X_2Y_2Z_2)$ is $1/2$. Assuming $v_2(n)=i>0$, we get the contribution
$$
\sum_{i=1\atop i\neq\alpha}^{\infty}\frac{d_k\left(2^{\text{min}(i,\alpha)}\right)d_l\left(2^i\right)d_m\left(2^{\text{min}(i,\alpha)}\right)}{2^{i+1}}
$$
when $i\neq \alpha$. When $i=\alpha$, we get the contribution
$$
\sum_{s=2}^{\infty}\frac{d_k\left(2^{\alpha+s}\right)d_l\left(2^{\alpha}\right)d_m\left(2^{\alpha+1}\right)}{2^{\alpha+s+1}}+
\sum_{t=2}^{\infty}\frac{d_k\left(2^{\alpha+1}\right)d_l\left(2^{\alpha}\right)d_m\left(2^{\alpha+t}\right)}{2^{\alpha+t+1}}.
$$
This gives the desired formula.

\subsection{Verifying $C_{k,l,m}(h)=\psi_{k,l,m}(h)$}
We recall that 
	\begin{align*}
		C_{k,l,m}(h)= \prod_{p\mid 2h}\left(\sum_{\nu_1,\nu_2,\nu_3\ge 0}d_{k-1}(p^{\nu_1})d_{l-1}(p^{\nu_2})d_{m-1}(p^{\nu_3})\frac{g(p^{\nu_1},p^{\nu_2},p^{\nu_3})}{[p^{\nu_1},p^{\nu_2},p^{\nu_3}]}\right),
	\end{align*}
where $g(p^{\nu_1},p^{\nu_2},p^{\nu_3})= 1$ if $\gcd (p^{\nu_1},p^{\nu_2}) \mid h, \gcd(p^{\nu_2},p^{\nu_3}) \mid h$ and $\gcd(p^{\nu_3},p^{\nu_1})\mid 2h$, otherwise $g(p^{\nu_1},p^{\nu_2},p^{\nu_3})= 0$. Let $p$ be an odd prime number and $v_p(h)=\alpha$. Then at most one of $\nu_1,\nu_2,\nu_3$ can be $>\alpha$. If
$\nu_1 >\alpha$, the contribution is
$$
\sum_{\nu_1>\alpha}\sum_{\nu_2=0}^{\alpha}\sum_{\nu_3=0}^{\alpha}\frac{d_{k-1}(p^{\nu_1})d_{l-1}(p^{\nu_2})d_{m-1}(p^{\nu_3})}{p^{\nu_1}}.
$$
Clearly $\sum_{\nu_2=0}^{\alpha} d_{l-1}(p^{\nu_2})=d_l(p^\alpha)$
and $\sum_{\nu_3=0}^{\alpha} d_{m-1}(p^{\nu_2})=d_m(p^\alpha)$.
Hence if $\nu_1 >\alpha$, applying \eqref{h}, we get that the contribution is
$$
\left(1-\frac{1}{p}\right)\sum_{i=1}^{\infty}\frac{d_{k}(p^{\alpha+i})d_l(p^{\alpha})d_m(p^{\alpha})}{p^{\alpha+i}}-\frac{d_k(p^{\alpha})d_l(p^{\alpha})d_m(p^{\alpha})}{p^{\alpha+1}}.
$$
Similarly, when $\nu_2 >\alpha$ and $\nu_3 >\alpha$, we get the contribution
$$
\left(1-\frac{1}{p}\right)\sum_{i=1}^{\infty}\frac{d_{k}(p^{\alpha})d_l(p^{\alpha+i})d_m(p^{\alpha})}{p^{\alpha+i}}-\frac{d_k(p^{\alpha})d_l(p^{\alpha})d_m(p^{\alpha})}{p^{\alpha+1}},
	$$
	and
	$$ 
	\left(1-\frac{1}{p}\right)\sum_{i=1}^{\infty}\frac{d_{k}(p^{\alpha})d_l(p^{\alpha})d_m(p^{\alpha+i})}{p^{\alpha+i}}-\frac{d_k(p^{\alpha})d_l(p^{\alpha})d_m(p^{\alpha})}{p^{\alpha+1}},
	$$
respectively. This leaves us with the case $\nu_1,\nu_2,\nu_3\le\alpha$.
Some careful computations in this case give the contribution
$$
\left(1-\frac{1}{p}\right)\sum_{i=0}^{\alpha-1}\frac{d_k(p^i)d_l(p^i)d_m(p^i)}{p^i}+\frac{d_k(p^{\alpha})d_l(p^{\alpha})d_m(p^{\alpha})}{p^{\alpha}}.
$$
Hence in $C_{k,l,m}(h)$, the Euler factor for an odd prime $p\mid h$ is exactly
$\psi_{k,l,m}(p^\alpha)$.

Now suppose $v_2(h)=\alpha$. The contribution when $0\le\nu_i\le\alpha$ for all $i=1,2,3$ is analogously 
	$$
	\frac{1}{2}\sum_{i=0}^{\alpha-1}\frac{d_k(2^i)d_l(2^i)d_m(2^i)}{2^i}+\frac{d_k(2^{\alpha})d_l(2^{\alpha})d_m(2^{\alpha})}{2^{\alpha}}.
	$$
Now we assume that $\nu_i \ge \alpha+1$ for at least one of $i=1,2,3$.
Now if $g(2^{\nu_1},2^{\nu_2},2^{\nu_3})=1$, then we have the following conditions: 
at most one of $\nu_1,\nu_2 >\alpha$, at most one of $\nu_2,\nu_3 >\alpha$ and at most one of $\nu_1,\nu_3 >\alpha+1$. Hence, as in \cite{MMS}, there are the following four possible choices:
	 \begin{enumerate}
	 	\item[i)] if $\nu_1=\alpha+1$, then $\nu_2 <\alpha+1$ and $\nu_3 \ge 0$,
	 	\item[ii)] if $\nu_1>\alpha+1$, then $\nu_2 <\alpha+1$ and $\nu_3 \le \alpha+1$,
	 	\item[iii)] if $\nu_1<\alpha+1$ and $\nu_2 \ge \alpha+1$, then $\nu_3 < \alpha+1$,
	 	\item[iv)] if  $\nu_1<\alpha+1$ and $\nu_2 < \alpha+1$, then $\nu_3 \ge \alpha+1$.
	 \end{enumerate}
In case i), the contribution is 
\begin{align*}
	 &\sum_{\nu_2\le \alpha}\sum_{\nu_3\le \alpha}\frac{d_{k-1}\left(2^{\alpha+1}\right)d_{l-1}\left(2^{\nu_2}\right)d_{m-1}\left(2^{\nu_3}\right)}{2^{\alpha+1}}+\sum_{\nu_2\le \alpha}\sum_{\nu_3\ge \alpha+1}\frac{d_{k-1}\left(2^{\alpha+1}\right)d_{l-1}\left(2^{\nu_2}\right)d_{m-1}\left(2^{\nu_3}\right)}{2^{\nu_3}}\\
	 =& \frac{d_{k-1}\left(2^{\alpha+1}\right)d_{l}\left(2^{\alpha}\right)d_{m}\left(2^{\alpha}\right)}{2^{\alpha+1}}+d_{k-1}\left(2^{\alpha+1}\right)d_{l}\left(2^{\alpha}\right)\sum_{i\ge \alpha+1}\frac{d_{m-1}\left(2^{i}\right)}{2^{i}}\\
	 =&\frac{1}{2}\sum_{i=1}^{\infty}\frac{d_{k-1}(2^{\alpha+1})d_l(2^{\alpha})d_m(2^{\alpha+i})}{2^{\alpha+i}},
	 \end{align*}
where the last step follows from $\eqref{h}$.
In cases ii), iii) and iv), similar computations lead to the contributions
$$
\frac{1}{2}\sum_{i=1}^{\infty}\frac{d_k(2^{\alpha+i})d_l(2^{\alpha})d_m(2^{\alpha+1})}{2^{\alpha+i}}-\frac{d_k(2^{\alpha+1})d_l(2^{\alpha})d_m(2^{\alpha+1})}{2^{\alpha+1}},
$$
$$
\frac{1}{2}\sum_{i=1}^{\infty}\frac{d_k(2^{\alpha})d_l(2^{\alpha+i})d_m(2^{\alpha})}{2^{\alpha+i}}-\frac{d_k(2^{\alpha})d_l(2^{\alpha})d_m(2^{\alpha})}{2^{\alpha+1}},
$$
and
$$
\frac{1}{2}\sum_{i= 1}^\infty \frac{d_k(2^{\alpha})d_l(2^{\alpha})d_m(2^{\alpha+i})}{2^{\alpha+i}}-\frac{d_k(2^{\alpha})d_l(2^{\alpha})d_m(2^{\alpha})}{2^{\alpha+1}},
$$
respectively. Note that the contributions in cases i) and iv) add up to
$$
\frac{1}{2}\sum_{i= 1}^\infty \frac{d_k(2^{\alpha+1})d_l(2^{\alpha})d_m(2^{\alpha+i})}{2^{\alpha+i}}-\frac{d_k(2^{\alpha})d_l(2^{\alpha})d_m(2^{\alpha})}{2^{\alpha+1}}.
$$
Adding all these, we get the desired formula.

\section*{Acknowledgements}
The authors thank Prof. M. Ram Murty for fruitful discussions and for his comments on an earlier version of this article. The authors also thank Prof. Maksym Radziwi\l\l  \ for some email exchanges regarding this work. The research of the second author was partially supported by IIT Delhi IRD project MI02455 and by SERB (ANRF) grant SRG/2022/001011.

\end{document}